\documentclass[psamsfonts]{amsart}

\newtheorem{theorem}{Theorem}[section]
\newtheorem{lemma}[theorem]{Lemma}
\newtheorem{proposition}[theorem]{Proposition}
\newtheorem{corollary}[theorem]{Corollary}

\theoremstyle{definition}
\newtheorem{definition}[theorem]{Definition}
\newtheorem{example}[theorem]{Example}

\newtheorem{remark}[theorem]{Remark}
\newtheorem{remarks}[theorem]{Remarks}

\theoremstyle{remark}

\numberwithin{equation}{section}

\newcommand{\NN}{\mathbb{N}}
\newcommand{\ZZ}{\mathbb{Z}}

\newcommand{\CC}{\mathbb{C}}

\newcommand{\PP}{\mathbb{P}}
\renewcommand{\AA}{\mathbb{A}}

\newcommand  {\shE}     {\mathcal{E}}
\newcommand  {\shF}     {\mathcal{F}}

\newcommand  {\shM}     {\mathcal{M}}

\newcommand  {\shL}     {\mathcal{L}}
\newcommand  {\shR}     {\mathcal{R}}

\newcommand  {\foa}     {\mathfrak{a}}

\newcommand  {\fod}     {\mathfrak{d}}

\newcommand  {\fom}     {\mathfrak{m}}

\newcommand  {\foq}     {\mathfrak{q}}

\newcommand  {\Char}    {\operatorname{char}}

\newcommand  {\Det}    {\operatorname{Det}}

\newcommand  {\dual}    {\vee}

\newcommand  {\Ext}     {\operatorname{Ext}}

\renewcommand  {\ker }  {\operatorname{kern}}

\newcommand  {\lra}     {\longrightarrow}

\renewcommand{\O}       {\mathcal{O}}

\newcommand  {\Proj}    {\operatorname{Proj}}

\newcommand  {\ra}      {\rightarrow}

\newcommand  {\Rel}     {\operatorname{Rel}}

\newcommand  {\Spec}    {\operatorname{Spec}}

\newcommand {\nonu} {\nu}

\newcommand {\leadno} { l}

\newcommand {\soclo} {\star}
\newcommand {\pasoclo} {\star}

\def\mydate{\number\day\space\ifcase\month \or January\or February\or March\or
April\or May\or
June\or July\or August\or September\or October\or November\or
December\fi \space\number\year}

\usepackage{amscd}
\usepackage{amssymb}

\begin{document}

\title[Tight closure and projective bundles]
{Tight closure and projective bundles}

% Remove or comment out any unused author tags.
% author one information

\author[Holger Brenner]{Holger Brenner}
\address{Mathematische Fakult\"at, Ruhr-Universit\"at, 
               44780 Bochum, Germany}
%\curraddr{}
\email{brenner@cobra.ruhr-uni-bochum.de}

%\thanks{}

\subjclass{14H60, 14J26, 13A35}

%\date{}

% at present the "communicated by" line appears only in ERA and PROC
%\commby{}

\dedicatory{\mydate}

\begin{abstract}
We formulate problems of tight closure theory in terms
of projective bundles and subbundles.
This provides a geometric interpretation of such problems and allows us
to apply intersection theory to them.
This yields new results concerning the tight closure of a
primary ideal in a two-dimensional graded domain.
\end{abstract}

\maketitle

%===========================================================
\section*{Introduction}

The aim of this paper is
to translate problems of tight closure theory
in terms of projective bundles and subbundles
in order to apply techniques of projective geometry such as intersection
theory to them.
This provides a geometric view on such problems
and enables us also to work often characteristic free.

We describe shortly the construction of the projective bundles arising
from tight closure.
The most basic problem of tight closure theory is to decide
whether $f_0 \in (f_1, \ldots ,f_n)^*$, where $f_0,f_1, \ldots ,f_n$
are elements in a Noetherian $K$-algebra $R$
($K$ a field of positive characteristic $p$).
This means by definition that there exists an element $c \in R$, not contained
in any minimal prime, such that
$cf_0^q \in (f_1^q, \ldots, f_n^q)$ for almost all powers $q=p^{e}$.

The starting point for our construction is the observation due to
Hochster (see \cite{hochstersolid}) that for a local complete
$K$-domain $(R,\fom)$ of dimension $d$ the
containment $f_0 \in (f_1, \ldots ,f_n)^*$
is equivalent to the property
that $H^d_\fom (A) \neq 0$ holds, where
$A=R[T_1, \ldots , T_n]/(f_1T_1+ \ldots +f_nT_n +f_0)$
is the so-called (generic) forcing algebra for the data
$f_1, \ldots, f_n;f_0$.
This characterization leads us to the study of the cohomological dimension
of the open subset
$D(\fom A) \subseteq \Spec \, A$;
if the ideal $(f_1, \ldots ,f_n)$ is primary to $\fom$, then this open
subset looks locally over $D(\fom) \subset \Spec \, R$
like an affine space $\AA^{d-1}$,
and the transition mappings are affine-linear.

In order to study the cohomological properties of this affine-linear bundle it
is helpful to embedd it into a projective
bundle.
This is achieved in the following way:
the spectrum
$\Spec\, R[T_0,T_1, \ldots , T_n]/(f_0T_0+f_1T_1+ \ldots +f_nT_n)$
yields a geometric vector bundle $V'$ over $D(\fom)$,
its sheaf of sections is given by the relations for the elements
$f_0,f_1, \ldots,f_n$.
The spectrum
$\Spec\, R[T_1, \ldots , T_n]/(f_1T_1+ \ldots +f_nT_n)$ yields a closed
subbundle $V \subset V'$ given by $T_0=0$.
These vector bundles yield
projective bundles $\PP(V) \subset \PP(V')$ and the complement
$\PP(V')-\PP(V)$ is isomorphic to our affine-linear bundle.

If $R$ is a graded normal domain
and if the $f_1, \ldots ,f_n$ are homogeneous $R_+$-primary
elements, we can go one step further and obtain a projective bundle
together with a projective subbundle of codimension one (called
forcing subbundle or forcing divisor) over $\Proj \, R$.
The cohomological dimension of the complement of the forcing divisor
is the same as the cohomological dimension of the affine-linear bundle
over $D(\fom)$, so we can work in an entirely projective
setting, which is moreover smooth whenever $R$ has an isolated singularity.

If $R$ is a normal standard graded domain of dimension two
then we are in a particularly manageable situation.
The construction leads to
projective bundles over smooth projective curves and
the question whether $f_0$ belongs to the tight closure of $(f_1, \ldots, f_n)$
is equivalent to the question whether the complement of the forcing divisor
is not an affine scheme (Proposition \ref{solid}).
This question is intimately related
to the question whether the forcing divisor is ample.

This geometric interpretation provides in particular a tool
to attack the following two problems of tight closure theory,
which we will encounter here for
several times and also in forthcoming papers
(\cite{tightelliptic},\cite{tightslope}).

The first problem is whether tight closure is the same as plus closure
in positive characteristic.
The plus closure of an ideal $I \subseteq R$ in a Noetherian domain
is just the contraction $I^+= R \cap IR^+$, where $R^+$ is the
integral closure of $R$ in an algebraic closure of $Q(R)$.
A positive answer to this problem would imply the localization problem
for tight closure.

If $R$ is graded, then the question whether
$f_0 \in (f_1, \ldots ,f_n)^{+{\rm gr}}$ 
is in our geometric setting equivalent to the existence
of projective subvarieties in $\PP(V')$ of $\dim R -1$ and disjoined
to the forcing divisor (\ref{plus}).
Therefore we look at the relation between tight closure and plus closure
as a relation
between intersection-geometric properties of the forcing divisor
and cohomological properties of the complement of it.

We will describe several situations in this paper where equality holds
(the results \ref{semiaffine}, \ref{pullbackample}, \ref{globalsection}
and \ref{numerischampel} are also true in characteristic zero,
whereas \ref{frobeniusartinschreier}, \ref{primarynumcrit}
and \ref{primarynumcor} need positive characteristic).
In \cite{tightelliptic} we will use our method to prove that
the tight closure and the plus closure
of a homogeneous $R_+$-primary ideal in a normal homogeneous
coordinate ring over an elliptic curve coincide in positive
characteristic.
The main ingredient for this result is the classification
of vector bundles on elliptic curves due to Atiyah,
which enables us to establish the same numerical criterion
for both properties.

The second problem is the graded Brian\c{c}on-Skoda-problem:
What is the minimal number $d_0$ such that
$R_{\geq d_0} \subseteq (f_1, \ldots, f_n)^*$ holds,
where $f_1, \ldots ,f_n$ are
homogeneous $R_+$-primary elements in a standard graded $K$-Algebra $R$?
It is known that this containment is true
for $d_0 \geq  d_1 + \ldots +d_n $, $d_i = \deg \, f_i$ and that this
number is a sharp bound in the parameter case,
see \cite{smithgraded} and \cite[Theorem 2.9 and 6.1]{hunekeparameter}.
However, this number is not much helpful in the general primary case.

Our interpretation suggests that in the two-dimensional situation the
number $(d_1 + \ldots + d_n)/(n-1)$ should be an important bound
for the degree,
since the top self intersection number of the forcing divisor is
$(d_1 + \ldots + d_n-(n-1)d_0) \deg\, \O_Y(1)$.
In this paper we show for $n=3$ (Theorem \ref{anwendung},
Theorem \ref{primarynumcrit}) that under some additional
conditions $(d_1+d_2+d_3)/2$ is the right bound.
This gives for example that $xyz \in (x^2,y^2,z^2)^*$
holds in $K[x,y,z]/(x^3+y^3+z^3)$.
In \cite{tightslope} we will show that
$R_{\geq (d_1+ \ldots +d_n)/(n-1)} \subseteq (f_1, \ldots, f_n)^*$
holds under the condition
that the relation bundle for $f_1, \ldots,f_n$ is
strongly semistable.
This rests upon
conditions for the inclusion in the tight closure 
in terms of the slopes of the corresponding
bundles.

The content of this paper is as follows.
In section \ref{forcingalgebra}
we discuss the characterization of tight closure
via solid closure
in terms of cohomological dimension and forcing algebras due to Hochster
\cite{hochstersolid}.
Henceforth we shall work rather with solid closure than with tight closure.
For two-dimensional rings this characterization leads to the
problem of affineness of open subsets (Proposition \ref{solidcd}).

The construction of the projective bundle, the forcing sequence
and the forcing divisor associated to a tight closure problem,
``is $f_0 \in (f_1, \dots ,f_n)^*?$'', and its basic properties
is given in section \ref{homogeneous}
and in section \ref{graded} for the graded case,
yielding bundles over $\Proj\, R$.

In section \ref{ample} we consider conditions for
the forcing divisor to be ample, to be basepoint free and to be big.
We give a geometric proof of
the result of Smith \cite[Theorem 2.2]{smithgraded}
that $f_0 \in (f_1, \ldots,f_n)^*$ implies $f_0 \in (f_1, \ldots ,f_n)$
for $\deg \, f_0 \leq \deg f_i,\, i=1, \ldots ,n$ and show
that it is also true
for solid closure (Corollary \ref{degreevanish}).

The rest of this paper is devoted to the study
of the tight closure of $R_+$-primary ideals
in a two-dimensional normal standard graded algebra $R$. We show that the
top self intersection number of the forcing divisor
is $(d_1 + \ldots +d_n - (n-1)d_0 ) \deg\, H$,
where $H$ is the hyperplane section on $\Proj \, R$, and that
this number is very important for the affineness of the complement
and hence for the tight closure question $f_0 \in (f_1, \dots ,f_n)^*?$

Sections \ref{ruled}-\ref{primary}
are concerned with the easiest case, the tight closure
of a homogeneous parameter ideal $(f_1,f_2)$.
Here our method brings rather new interpretations and proofs than
new results.
The construction yields ruled surfaces over the
corresponding smooth projective curve together with a forcing section
(Corollary \ref{forcingrule}).
The tight closure problem becomes a question on the ampleness
of this divisor (Theorem \ref{ampelcrit})and the number $(d_1+d_2-d_0) \deg \, H$
is the self intersection number of it.
We recover the so-called vanishing theorem that
$(f_1,f_2)^*=(f_1,f_2) + R_{\geq d_1 +d_2}$ holds for $p=0$ or $p >>0$
(Corollary \ref{degreecrit2}).

In giving examples of ruled surfaces arising from forcing data
we encounter
Hirzebruch surfaces, incidence varieties,
a classical construction of Serre of a Stein but non-affine variety
and a new class of counterexamples to the hypersection problem
(sections \ref{examples} and \ref{complex}).
This shows also that we can establish
geometrically interesting properties from results of tight closure theory. 

Section \ref{positive} deals with the plus closure and section
\ref{primary} with primary relations and how they influence the
$e$-invariant of the ruled surfaces.

Section \ref{higherrank} then deals with the tight closure of three primary
homogeneous elements $f_1,f_2,f_3$ in a two-dimensional graded ring,
yielding projective bundles of rank two over the curve.
This is already a very subtile situation where new phenomena occurs,
and a detailed study of the geometric situation is necessary to obtain
results on tight closure.
If the number $d_1 +d_2 + d_3 - 2d_0 $ is $ \leq 0$,
then under some extra conditions on the $f_1,f_2,f_3$ concerning
their relations, we show
that $R_{\geq d_0} \subseteq (f_1,f_2,f_3)^*$
(Theorem \ref{anwendung} and also in the plus closure,
Theorem \ref{primarynumcrit}).
To mention just one example, it follows that for $R=K[x,y,z]/(x^d+y^d+z^d)$
we obtain
$R_{\geq d} \subseteq (x^{d_1},y^{d_2},z^{d_3})^*$
for $d_1+d_2+d_3 =2d$, $d_i < d$ (Example \ref{fermatdegree}).

I thank Manuel Blickle for many suggestions
and the referee for careful reading and useful and critical remarks which
helped substantially to improve this paper.

%===========================================================
\section{Forcing algebras and cohomological dimension}
\label{forcingalgebra}

Let $R$ denote a commutative ring and let $f_1, \ldots ,f_n \in R$ and
$f_0 \in R$ be elements. The $R$-algebra
$$A=R[T_1,\ldots ,T_n]/(f_1T_1+\ldots +f_nT_n+f_0) \, ,$$
is called the {\em {\rm(}generic{\rm)} forcing algebra}
for the elements $f_1, \ldots ,f_n;f_0$.
The forcing algebra forces that
$f_0 \in (f_1, \ldots, f_n)A$ and every other $R$-algebra with this property
factors through $A$.
For studying tight closure problems in terms of forcing algebras,
the cohomological properties of 
the subsets $D(\fom A) \subseteq \Spec \, A$ for the maximal
ideals $\fom \in \Spec\, R$ are important.
Recall that the {\em cohomological dimension} $cd \,(U)$
of a scheme $U$ is the maximal number $i$ such
that there exists a quasicoherent sheaf $\shF$ with $H^{i}(U, \shF) \neq 0$
(see \cite{hartshornecohdim} for this notion).
For an ideal $\foa \subseteq R$ we call the maximal number $j$
such that there exists an $R$-module $M$ with
$H^{j}_\foa (M) \neq 0$ the {\em cohomological height}, $ch\, (\foa)$
(this is sometimes called the local cohomological dimension).
For $\dim \, R \geq 2$ we have $ch \,(\foa) = cd \,(D(\foa)) +1 $,
due to the long exact sequence of local cohomology.

We will not go back to the definition of tight closure but
we recall the notion of {\em solid closure} in a form which is suitable
for our purpose. 

\begin{definition}
Let $R$ be a Noetherian ring and let $f_1, \ldots ,f_n, f_0 \in R$.
Then $f_0$ belongs to the {\em solid closure},
$f_0 \in (f_1, \ldots ,f_n)^{\soclo}$, if and only if
for every local complete domain $R' =  \hat{R_{\fom}}/\foq $
(where $\fom $ is a maximal ideal of $R$ and $\foq$ is a minimal prime
of $\hat{R_{\fom}}$) we have
$H^d_{\fom'}(A') \neq 0$, where $d = \dim \, R'$ and $A'$ is the
forcing algebra over $R'$.
\end{definition}

\begin{remarks}
\label{solidremarks}
This definition coincides with the definition given in
\cite[1.2]{hochstersolid} due to
\cite[Corollary 2.4 and Proposition 5.3]{hochstersolid}.
The condition must only be checked for the maximal ideals
$\fom \supseteq (f_1, \ldots ,f_n)$.

Suppose that $R$ contains a field of characteristic $p >0$
and suppose furthermore that $R$
is essentially of finite type over an excellent local ring
or that the Frobenius  endomorphism is finite.
Then the tight closure of an ideal is the same as its solid closure,
see \cite[Theorem 8.6]{hochstersolid}
(this is not true in characteristic $0$ for $\dim \, R \geq 3$).
\end{remarks}

\begin{proposition}
\label{solidcd}
Let $R $ be a normal excellent domain.
Let $f_1, \ldots ,f_n \in R$ be primary to a maximal ideal $\fom$
of height $d$,
$f_0 \in R$ and let
$A=R[T_1,\ldots ,T_n]/(f_1T_1+\ldots +f_nT_n+f_0)$ be the forcing algebra.
Then the following hold.

\renewcommand{\labelenumi}{(\roman{enumi})}
\begin{enumerate}

\item
$f_0 \in (f_1, \ldots, f_n)^\soclo$ if and only if the
cohomological height
of the extended ideal $\fom A$ is $d$.

\item
If $d \geq 2$, then $f_0 \in (f_1, \ldots, f_n)^\soclo$ if and only if
the cohomological dimension of $W=D(\fom A) \subset \Spec\, A$
is $d-1$.

\item
If $d=2$,
then $f_0 \in (f_1,\ldots, f_n)^\soclo $
if and only if $\, D(\fom A)$ is not an affine scheme.

\end{enumerate}
\end{proposition}
\proof
(i)
Since the completion of a normal and excellent domain is again a domain
the condition $f_0 \in (f_1, \ldots ,f_n)^\soclo$
is equivalent to $H^d_{\fom '}(A') \neq 0$, where $R'$ is the completion of
$R_{\fom}$. Since cohomology commutes under completion
this is equivalent to $H^d_\fom (A) \neq 0$.
Since $H^d_\fom (A) = H^d_{\fom A} (A)$
this implies that $ch \, (\fom A) \geq d$, and equality must hold
since the cohomological height of $\fom A$ can not be bigger
than $ch \, (\fom) = d$.
On the other hand, if $H^d_\fom (A) =0$, then this holds for every
$A$-module $M$, since $A^{(J)} \ra M \ra 0$ and since $H^{d+1}_\fom(-)=0$.

(ii) follows from (i) by the long exact sequence of local cohomology.
(iii) follows from (ii) and
the cohomological characterization of affine schemes.
\qed

\begin{lemma}
\label{finiteextension}
Let $R$ be a Noetherian ring, 
let $f_1, \ldots ,f_n,f_0 \in R$ be elements and set
$A=R[T_1,\ldots ,T_n]/(f_1T_1+\ldots +f_nT_n+f_0)$.
Let $W=D(\fom A) \subset \Spec\, A$.
Then the following hold.
\renewcommand{\labelenumi}{(\roman{enumi})}
\begin{enumerate}

\item
Suppose that $R, \fom$ is local of
dimension $\dim \, R=d \geq 2$
and suppose that
there exists another local Noetherian ring $R'$ of dimension $d$
and a ring homomorphism $R \ra R'$ such that
$V(\fom R') =V(\fom_{R'})$ and
$f_0 \in (f_1,\ldots,f_n)R'$ hold.
Then the cohomological dimension of $\, W$ is $d-1$
{\rm (}and $f_0 \in (f_1, \ldots ,f_n)^\soclo$, if $R$ is normal and
excellent{\rm )}.

\item
Let $R$ be a normal domain over a field  $K$ of characteristic zero.
If there exists a finite extension $R \subseteq R'$ such that
$f_0 \in (f_1, \ldots, f_n)R'$, then
already $f_0 \in (f_1, \ldots, f_n)R$.
\end{enumerate}
\end{lemma}

\proof
(i)
The morphism $\Spec \, R' \ra \Spec \, R$ lifts to a
morphism $\varphi: \Spec\, R' \ra \Spec\, A$ and
$\varphi^{-1} (W) = D(\fom_{R'})$.
Thus we have an affine morphism $D(\fom_{R'}) \ra W$
and the cohomolgical dimension of $D(\fom_{R'})$ is
$d-1$.

(ii)
This follows from the existence of the trace map, see
\cite[Remarks 9.2.4]{brunsherzog}.
\qed

\begin{remark}
In positive characteristic it is sometimes possible to
show that $f_0 \in (f_1, \ldots ,f_n)^*$ by giving a finite
extension $R \subseteq R'$ where $f_0 \in (f_1 , \ldots ,f_n)R'$ holds.
In fact it is a tantalizing question of Hochster
whether this is always true -- i. e. whether tight closure is the same
as plus closure.
A result of Smith \cite{smithparameter} says that this is
true for parameter ideals.
In \cite{tightelliptic} we show that this is also true for homogeneous
$R_+$-primary ideals in an affine normal cone over an elliptic curve.
\end{remark}

The {\em superheight} of an ideal $\foa \subseteq R$
is the maximal height of $\foa R'$ in any Noetherian $R$-algebra $R'$.
The superheight of an ideal is less or equal its cohomological
height, and solid closure gives in characteristic $0$
examples that it may be less.
In particular tight closure gives examples of open subsets
$D(\foa)$ such that $\foa$ has superheight one, but $D(\foa)$
is not affine.
For other examples see \cite{neeman} and \cite{brennersuperheight}.
We will apply this in section \ref{complex} to give new counterexamples
to the hypersection problem of complex analysis.

\begin{corollary}
\label{superhoehe}
Let $K$ be a field of characteristic zero and let
$R$ be a normal excellent $K$-domain.
Let $f_1,\ldots, f_n$ be primary to a maximal ideal $\fom$
of height $d$
and let $f_0 \in R$.
Suppose that
$f_0 \not\in (f_1, \ldots, f_n)$, but $f_0 \in (f_1, \ldots, f_n)^\soclo$.
Then the cohomological height of
$\fom A \subseteq A= R[T_1, \ldots ,T_n]/ (\sum f_iT_i+f_0)$
is $d$ and its superheight is $< d$.
\end{corollary}

\proof
Let $R'$ denote a local normal Noetherian domain
of dimension $n \leq d$ and let
$ \varphi: A \ra R'$
be a homomorphism such that
$ V( (\fom_R A)R')= V(\fom_{R'})$.
This gives a homomorphism
$\psi: R \ra R'$ such that
$V(\fom_R R') = V(\fom_{R'})$.
If $n = d$, then $\hat{R_\fom} \ra \hat {R'}$ would be finite
(after enlarging the base field) and \ref{finiteextension}(ii)
would give $f_0 \in (f_1, \ldots, f_n)$.
Hence the superheight is $< d$, but the cohomological height
is $d$ due to \ref{solidcd} (i).
\qed

\section{Homogeneous forcing algebras and projective bundles}
\label{homogeneous}

Consider the mapping
$\Spec \,A \longrightarrow \Spec \, R$ over the open
subset $U=D(I)=D(f_1,\ldots ,f_n)$, where $A$ is the forcing algebra
for $f_1, \ldots ,f_n; f_0$.
On $D(f_i)$, $i \geq 1$, one can identify
$$ A_{f_i} = (R[T_1,\ldots ,T_n]/(f_1T_1+\ldots +f_nT_n +f_0))_{f_i} \cong
R_{f_i}[T_1, \ldots ,T_{i-1},T_{i+1},\ldots ,T_n] \, .$$
So this mapping looks locally like
$D(f_i) \times \AA ^{n-1} \longrightarrow D(f_i)$.
The transition mapping on $D(f_if_j)$ is given by
$$R_{f_if_j}[T_1,\ldots ,T_{i-1},T_{i+1},\ldots ,T_n] \longrightarrow
R_{f_if_j}[T_1,\ldots ,T_{j-1},T_{j+1},\ldots ,T_n] $$
where $ T_k \longmapsto T_k $ for $k \neq i,j $ and
$T_j \longmapsto -1/f_j(\sum_{i \neq j}f_iT_i +f_0) $.
This is an affine-linear mapping.
Therefore we say that the forcing bundle $\Spec \, A|_{D(I)}$
is an {\em affine-linear bundle} of rank $n-1$.
It is not a vector bundle in general.

We show how to associate to elements $f_1, \ldots ,f_n;f_0$
a projective bundle over $D(I)$
together with a projective subbundle of codimension one such that
the complement of the subbundle is the affine-linear bundle.
This is more generally possible for every affine-linear bundle.

\begin{proposition}
\label{forcingsequence1}
Let $R$ be a commutative ring and let
$f_1,\ldots, f_n$ and $f_0$ be elements and set
$I=(f_1,\ldots,f_n)$, $U=D(I)$.
The schemes
$$ V= \Spec R[T_1,\ldots, T_n]/(\sum_{i=1}^n f_iT_i)|_U
\, \mbox{ and } \,
V' = \Spec R[T_0,\ldots, T_n]/(\sum_{i=0}^n f_iT_i)|_U $$
are vector bundles on $U$.
They are related by the short exact sequence of vector bundles
$$0 \longrightarrow V \longrightarrow V' 
\stackrel{T_0}{\longrightarrow} \AA^1_U \longrightarrow 0 \, .$$
The inclusion $V \subset V'$ yields
a closed embedding
$\PP(V) \hookrightarrow \PP(V')$ of projective bundles over $U$.
Its complement
$\PP(V') -\PP(V)$ is isomorphic to the
forcing affine-linear bundle
$$\Spec\, R[T_1,\ldots,T_n]/(f_1T_1+\ldots +f_nT_n+f_0)|_U \, .$$
\end{proposition}

\proof
The bundle $V$ is on $D(f_i)$, $i=1, \ldots ,n$,
isomorphic to
$$\Spec R_{f_i}[T_1, \ldots,T_{i-1},T_{i+1}, \ldots ,T_n] \, ,$$
and the transition functions send
$T_i \mapsto
-1/f_i(f_1T_1+ \ldots +f_{i-1}T_{i-1}+f_{i+1}T_{i+1} + \ldots +f_nT_n)$,
thus they are linear and $V$ (and $V'$) is a vector bundle
on $U=\bigcup_{i=1}^n D(f_i)$.
 
The linear form $T_0$ is a global function on $V'$ which yields a
linear mapping to $\AA^1_U$. Its zero set is $V$.
Looking at $D(f_i)$, the exactness of the sequence is clear.

$\PP(V')$ is the projective bundle corresponding to the geometric
vector bundle $V'$.
The cone mapping $V' \dashrightarrow \PP(V')$ maps $V(T_0-1)$
isomorphically onto $\PP(V')-\PP(V)$.
\qed

\begin{definition}
We call the short exact sequence in \ref{forcingsequence1} the
{\em forcing sequence} and we call $\PP(V')$ the {\em projective bundle}
and $\PP(V)$ the {\em forcing projective subbundle}
or the {\em forcing divisor} associated to the elements $f_1, \ldots , f_n;f_0$.
\end{definition}

\begin{remark}
The sections $\Spec \, R \ra \Spec\, R[T_1, \ldots ,T_n]/(\sum_{i=1}^n f_iT_i)$
are the relations for the ring elements $f_1, \ldots ,f_n$.
This is true for every open subset in $\Spec \, R$.
We call this sheaf of sections the {\em sheaf of relations}
$\shR=\Rel(f_1,\ldots,f_n)^{\tilde{ }} \,$.
On $U=D(I)$, this is a locally free sheaf and we get
the short exact forcing sequence of locally free sheaves
$$0 \longrightarrow \shR \longrightarrow \shR' 
\longrightarrow \O_U \ra 0 \, .$$
These extensions are classified by $H^1(U,\shR)= \Ext^1(\O_U,\shR)$.
The elements $f_0$ and $f_1, \ldots ,f_n$ define the
$\check{\rm C}$ech-cocycle
$$(0, \ldots, -f_0/f_i, 0, \ldots ,0,f_0/f_j,0, \ldots, 0)
\in \Gamma(D(f_if_j), \shR) \subseteq \Gamma(D(f_if_j), \O_U^n) \, .$$
The dual sheaf $\shF =\shR^\dual$ is the {\em sheaf of linear forms}
for the vector bundle $V$, thus $V= \Spec \, S(\shF)$
and $\PP(V)= \Proj \, S(\shF)$.
Geometric vector bundles, their sheaf of relations and their sheaf of
linear forms are essentially equivalent objects; in this paper
we shall take mostly the viewpoint of geometric vector bundles,
since in this form they appear starting from forcing algebras.
\end{remark}

We gather together some characterizations of
$f_0 \in (f_1, \ldots ,f_n)$ in terms of
the geometric objects we consider.

\begin{lemma}
\label{trivial}
Let $R$ be a commutative ring and let
$f_0,f_1, \ldots,f_n \in R$ be elements,
$U=D(f_1,\ldots,f_n)$.
Let $A=R[T_1, \ldots, T_n]/(f_1T_1+ \ldots +f_nT_n+f_0)$
be the forcing algebra. Then the following are equivalent.

\renewcommand{\labelenumi}{(\roman{enumi})}
\begin{enumerate}
\item 
$f_0 \in (f_1,\ldots,f_n)$.

\item
There exists a section $\Spec \, R \ra \Spec \, A$.

\item
The forcing algebra $A$ is isomorphic to
the algebra of relations
$$R[T_1, \ldots ,T_n]/ (f_1T_1+ \ldots +f_nT_n) \, .$$
Suppose furtheron that $R =\Gamma(U,\O_X)$ {\rm(}e.g. if $R$ is normal
and ${\rm ht}\, I \geq 2 $.{\rm )}
Then these statements are also equivalent with

\item
The affine-linear bundle $\Spec \, A|_U$ has a section over $U$.

\item
There exists a section $U \ra \PP(V')$ which does not meet
$\PP(V)$.

\item
The forcing sequence splits.

\item
The elements $f_1, \ldots ,f_n,f_0$ define the zero element in $H^1(U, \shR)$.

\end{enumerate}
\end{lemma}

\proof
Suppose (i) holds, say $-f_0=\sum a_if_i$.
Then $T_i \mapsto T_i+a_i$ is well defined and gives the isomorphism in (iii).
On the other hand, a relation algebra has the zero section, thus the first
three statements are equivalent.

(ii) $\Rightarrow$ (iv) is a restriction, and (iv) $ \Rightarrow $ (ii)
is true under the additional assumption. (iv) and (v) are equivalent
due to \ref{forcingsequence1}.

(i) gives also directly a section for $V' \ra \AA_U \ra 0$, thus we
get (vi), which is equivalent with (vii).
If the sequence splits, then $V' = V \oplus \AA$ on $U$ and the complement
of $\PP(V)$ is the vector bundle $V$, which has the zero section.
\qed

\section{The graded case: bundles on projective varieties}
\label{graded}

In order to use methods of projective geometry such as intersection theory
to study the affineness
of an open subset inside the spectra of a forcing algebra,
we stick now to the graded case, where we get projective bundles
over projective varieties.

Let $K$ be a field and let $R$ be a
standard $\NN$-graded $K$-algebra, i.e. $R_0=K$ and $R$ is generated
by finitely many elements of first degree.
Let $f_i$ be homogeneous elements of $R$ of degrees $d_i$.
We say that the $f_i$ are primary if $D(R_+) \subseteq D(f_1, \ldots ,f_n)$.
We may find degrees $e_i$ (possibly negative) for $T_i$ such that
the polynomials
$\, \sum_{i=1}^nf_iT_i \,$, $\, \sum_{i=0}^n f_iT_i \,$ and
$\, \sum_{i=1}^n f_iT_i+f_0$ are homogeneous
(for the last polynomial $e_i=d_0-d_i$ is the only choice).

Let
$A=R[(T_0),T_1, \ldots ,T_n]/(P)$, where $P$ is one of these polynomials.
Then $A$ is also graded and we have the following commutative diagram.

$$
\begin{CD}
\Spec \, A \supset D(R_+A) \, \,   @>  >> \, \,  D_+(R_+ A) \subset \Proj\, A   \\
\, \,  @VVV     @VVV \\
\Spec \, R  \supset D(R_+) \, \,  @>  >> \Proj \, R
\end{CD}
$$

\medskip
\noindent
For $Y= \Proj\, R$ and a number $m$ we set
$$\AA_Y(m) := D_+(R_+) \subset \Proj\, R[T] ,\, \,  \deg \, T =  -m \, .$$
This line bundle on $Y$ is also $\AA_Y(m) = \Spec \, S(\O_Y(m))$
and its sheaf of sections is $\O_Y(-m)$.
Thus $\AA_Y(1)$ is the tautological bundle.
(The algebras $A$ may have negative degrees,
but $\Proj \, A$ can be defined as well,
see \cite{brennerschroer}. The open subset
$D_+(R_+A) \subset \Proj \, A$ is the same as
$D_+(R_+A) \subseteq \Proj \, A_{\geq 0}$.)

\begin{proposition}
\label{buendelalsproj}
Let $R$ be a standard graded $K$-algebra and let
$f_1,\ldots,f_n$ be homogeneous primary elements.
Let $d_i = \deg \, f_i $ and
fix a number $m \in \ZZ$ and set $e_i =m-d_i$.
Then the following hold.

\renewcommand{\labelenumi}{(\roman{enumi})}
\begin{enumerate}

\item
Set $\deg \, T_i =e_i$. Then
$$ \Proj\, R[T_1,\ldots,T_n]/(\sum_{i=1}^n f_iT_i) \, \, \supset \, \, D_+(R_+)
\longrightarrow \Proj\, R$$
is a vector bundle $V_m$ of rank $n-1$ over $Y=\Proj\, R$.

\item
For this vector bundle $V_m$ we have the exact sequence
of vector bundles
$$0 \longrightarrow V_m \longrightarrow
\AA_Y(-e_1) \times_Y \ldots \times_Y \AA_Y(-e_n)
\stackrel{\sum f_i}{\lra} \AA_Y(-m) \longrightarrow 0 $$

\item
We have $\Det \, V_m \cong \AA_Y(-\sum_{i=1}^n e_i +m)
= \AA_Y(\sum_{i=1}^n d_i -(n-1)m) $.

\item
We have $V_{m'}= V_m \otimes \AA_Y(m-m')$.

\item
The projective bundle
$\PP(V_m)$ does not depend on the chosen degree $m$.
For the relatively very ample sheaf
$\O_{\PP(V_m)}(1) $ on $\PP(V_m)$ we have
$$j^* \O_{\PP(V_{m'})} (1) =  \O_{\PP(V_m)} (1) \otimes \pi^*\O_Y(m-m')\, ,$$
where
$j: \PP(V_m) \ra \PP(V_m \otimes  \AA_Y(m-m'))$
is the isomorphism and $\pi: \PP(V_m) \ra Y$ is the projection.

\end{enumerate}
\end{proposition}

\proof
(i) and (ii).
First note that the natural mapping ($\deg\, T_i=e_i$)
$$\Proj \, R[T_1, \ldots ,T_n] \supseteq D_+(R_+) \lra
\AA_Y(-e_1) \times_Y \ldots  \times_Y \AA_Y(-e_n)$$
is an isomorphism.
The ring homomorphism
$R[T] \ra R[T_1,\ldots, T_n], \, T\mapsto \sum f_iT_i$
is homogeneous for $\deg \, T =m$.
This gives the epimorphism of vector bundles, since the $D_+(f_i)$
cover $Y$.
Its kernel is given by
$ D_+(R_+) \subset  \Proj \,  R[T_1,\ldots, T_n]/(\sum f_iT_i)$,
thus this is also a vector bundle on $Y$.

(iii) follows from (ii).
If we
tensorize the exact sequence for $V_m$ with $\AA_Y(m-m')$ we
get the sequence for $V_{m'}$, hence (iv) follows.

(v)
$\PP(V)$ does not change when $V$ is tensorized
with a line bundle. The relatively very ample sheaves behave like stated due
to \cite[Proposition 4.1.4]{EGAII}.
\qed

\begin{remark}
We denote by $\shR(m)$
the locally free sheaf of sections in the vector bundle $V_m$.
This is the sheaf of relations of total degree $m$,
and $\shR(m) \otimes \O(m'-m)=\shR(m')$ holds.
The sequence in \ref{buendelalsproj} yields
the short exact sequence
$$ 0 \lra \shR (m)  \lra \oplus_i \O_Y(e_i)
\stackrel{\sum f_i}{\lra} \O_Y(m) \lra 0 \, .$$
The sheaf of linear forms of total degree $m$
is the dual sheaf $\shF(-m)= \shR(m)^\dual$,
thus $V_m= \Spec \, S(\shF(-m))$ and $\PP(V)=\PP(\shF)$.
The corresponding sequence is
$$0 \lra \O_Y(-m) \stackrel{f_1, \ldots ,f_n}{\lra}
\oplus_i \O_Y(-e_i) \lra \shF (-m) \lra 0 \, .$$
The most important choice for $m$ will be $m=d_0$, where $d_0$
is the degree of another homogeneous element $f_0$.
\end{remark}

\begin{remark}
\label{chern}
The sequence in \ref{buendelalsproj} (ii) allows us
to compute inductively the
Chern classes of the vector bundles $V_m$ (or its sheaf of linear forms
$\shF(-m)$).
For the Chern polynomial $c_t(V_m)= \sum_{i}c_i(V_m)t^{i}$
we get the relation
(let $H$ denote the hyperplane section of $Y$)
$$c_t(V_m)(1-mHt)=(1-e_1Ht) \cdots (1-e_nHt) \, .$$
This yields
$c_0(V_m)=1,\, c_1(V_m)= (-e_1- \ldots -e_n+m)H,\,
c_2(V_m)=(\sum_{i_1,i_2} e_{i_1}e_{i_2}-(e_1+ \ldots +e_n-m)m) H.H $ etc.
\end{remark}

\begin{proposition}
\label{forcingsequence2}
Let $R$ be a standard graded $K$-algebra, let
$f_1,\ldots,f_n$ be homogeneous primary elements and let
$f_0 \in R$ be also homogeneous.
Let $d_i = \deg \, f_i $ and fix a number $m \in \ZZ$.
Let $\deg \, T_i =e_i =m-d_i$.
Let
$$V_m=D_+(R_+) \subset \Proj\, R[T_1,\ldots,T_n]/(\sum_{i=1}^n f_iT_i) 
\, \, \mbox{ and }$$
$$V'_m= D_+(R_+) \subset \Proj\, R[T_0,\ldots,T_n]/(\sum_{i=0}^n f_iT_i)$$
be the vector bundles on $Y=\Proj\, R$ due to {\rm \ref{buendelalsproj}}.
Then the following hold.

\renewcommand{\labelenumi}{(\roman{enumi})}
\begin{enumerate}

\item 
There is an exact sequence of vector bundles on $Y$,
$$0 \longrightarrow V_m \longrightarrow V'_m
\stackrel{T_0}{\longrightarrow} \AA_Y(-e_0) \longrightarrow 0 \, .$$

\item
The embedding $\PP(V) \hookrightarrow \PP(V')$ does not depend
on the degree $m$ {\rm (}and we skip the index $m$ inside $\PP(V)${\rm )}.
The complement of $\PP(V)$ is
$$\PP(V') -\PP(V) \cong D_+(R_+)
\subseteq \Proj \, R[T_1, \ldots ,T_n]/(f_1T_1+ \ldots +f_nT_n+f_0) \, .$$

\item
Let $E$ be the Weil divisor {\rm (}the hyperplane section{\rm)} on $\PP(V')$
corresponding to the relatively very ample
invertible sheaf $\O_{\PP(V')} (1)$ {\rm (}depending of the degree{\rm )}.
Then we have the linear equivalence of divisors
$\PP(V) \sim E +e_0 \pi^*H$,
where $H$ is the hyperplane section of $\, Y$.
If $e_0=0$, then $\PP(V)$ is a hyperplane section.

\item
The normal bundle for
$\PP(V) \hookrightarrow \PP(V')$ on $\PP(V)$ is
$\AA_{\PP(V)}(-1) \otimes \pi^*\AA_Y(-e_0)$.

\end{enumerate}
\end{proposition}
\proof
(i).
The homogeneous ring homomorphisms
$$R[T_0] \, \lra \, R[T_0, \ldots, T_n]/(\sum_{i=0}^n f_iT_i)
\, \, \, (\deg \,T_0= e_0)  \, \mbox{ and }$$
$$ R[T_0, \ldots, T_n]/(\sum_{i=0}^n f_iT_i) \, \lra \,
R[T_1, \ldots, T_n]/(\sum_{i=1}^n f_iT_i), \, \, \, T_0 \longmapsto 0 $$
induce the morphisms on $D_+(R_+)$.
The exactness is clear on $D(f_i)$, $i=1, \ldots ,n$, and they cover
$D_+(R_+)$.

The first statement in (ii) is clear, thus we assume $e_0=0$.
The homogeneous ring homomorphism
$R[T_0, \ldots, T_n]/(\sum_{i=0}^n f_iT_i) \ra
R[T_0, \ldots, T_n]/(\sum_{i=1}^n f_iT_i +f_0)$ where $T_0 \mapsto 1$
yields the closed embedding
$$ \Proj \,  R[T_0, \ldots, T_n]/(\sum_{i=1}^n f_iT_i +f_0) \supseteq
D_+(R_+) \hookrightarrow V'_m \, , $$
where the image is given by $T_0 =1$.
But this closed subset $V_+(T_0-1) \subseteq V'_m$
is isomorphic to $\PP(V') -\PP(V)$
under the cone mapping $V'_m \dashrightarrow \PP(V')$.

(iii).
The mapping
$T_0: V'_m \ra \AA_Y(-e_0)$ yields via the tautological
morphism $\AA_{\PP(V'_m)}(1) \ra V'_m$ a morphism of line bundles on
$\PP(V'_m)$,
$\AA_{\PP(V'_m)}(1) \ra \pi^*\AA_Y(-e_0)$.
This corresponds to a section in the line bundle
$\AA_{\PP(V'_m)}(-1) \otimes \pi^*\AA_Y(-e_0)$
with zero set $\PP(V_m)$.
Thus $\PP(V_m) \sim E +e_0 \pi^*H$.

(iv).
Let $i:\PP(V) \hookrightarrow \PP(V')$ be the inclusion.
Then
$i^*(\AA_{\PP(V')}(-1) \otimes \pi^*\AA_Y(-e_0))
=\AA_{\PP(V)}(-1) \otimes q^*\AA_Y(-e_0)$
is the normal bundle on $\PP(V)$.
\qed

\begin{definition}
\label{forcingdef}
We call the sequence in \ref{forcingsequence2} (i) again
the {\em forcing sequence} and we
denote the situation $\PP(V) \hookrightarrow \PP(V')$
by $\PP(f_1, \ldots, f_n; f_0)$.
This is a projective bundle of rank $n-1$ together
with the forcing divisor $\PP(V)= \PP(f_1, \ldots ,f_n)$
over $Y$.
\end{definition}

\begin{remark}
\label{forcingsheaf2}
Corresponding to the forcing sequence of vector bundles in
\ref{forcingsequence2}
we have the exact sequence of relations
$0 \ra \shR(m) \ra \shR'(m) \ra \O_Y(e_0) \ra 0 $
of total degree $m$.
For $e_0=0$ (or $m=d_0$) this extension corresponds to a cohomology
class $c \in H^1(Y, \shR(m))$.
The forcing sequence for the linear forms is
$ 0 \ra \O_Y(-e_0)  \ra \shF'(-m) \ra \shF (-m)  \ra 0 $.
\end{remark}

The next results show that we can express the properties which are of interest
from the tight closure point of view in terms of the projective bundles
on $Y$.

\begin{lemma}
\label{trivialtwo}
In the situation of {\rm \ref{forcingsequence2}} the following
are equivalent.

\renewcommand{\labelenumi}{(\roman{enumi})}
\begin{enumerate}

\item 
$f_0 \in (f_1, \ldots ,f_n)$.

\item
There is a section $Y \ra \PP(V')$ disjoined to $\PP(V) \subset \PP(V')$.

\item
The forcing sequence
$0 \ra V_m \ra V'_m \ra \AA_Y(-e_0) \ra 0 $ splits.

\item
Let $e_0=0$.
The corresponding cohomological class in $H^1(Y,\shR (m))$ vanishes.
\end{enumerate}
\end{lemma}

\proof
Suppose that (i) holds and write $-f_0 =\sum_{i=1}^n a_if_i$,
where the $a_i$ are homogeneous. Set $e_0 =0$.
The $a_i$ define a homogeneous mapping
$$R[T_0,T_1,\ldots,T_n]/(\sum_{i=0}^n f_iT_i) \lra R
\mbox{ by }
T_0 \ra 1,\, \,  T_i \ra a_i \, .$$
The corresponding mapping
$Y \ra V'_m$ induces $Y \ra \PP(V')$ and its image
is disjoint to $\PP(V)$.

Suppose that (ii) holds.
A section in $\PP(V'_m)$ corresponds to a line bundle $L$ on $Y$
and an embedding $L \hookrightarrow V'_m$,
see \cite[Proposition 7.12]{haralg}.
Since the section is disjoined to $\PP(V)$, the
morphism $V_m \oplus L \ra V'_m$ is an isomorphism,
hence the sequence splits.

(iii).
The splitting yields a section
$\AA_Y(-e_0) \ra V'_m$ and this means a homogeneous mapping
$ R[T_0,T_1,\ldots,T_n]/(\sum_{i=0}^n f_iT_i) \ra R[T_0]$.
For $T_0=1$ we get a solution for (i).
(iii) and (iv) are equivalent. 
\qed

\begin{example}
\label{null}
Let $R$ denote a standard graded $K$-algebra
and let $f_1, \ldots ,f_n $ be homogeneous primary elements of degrees $d_i$.
Let $f_0=0$.
Then
$$ R[T_0, \ldots , T_n]/(\sum_{i=0}^n f_iT_i)
= R[T_1, \ldots , T_n]/(\sum_{i=1}^n f_iT_i) [T_0]$$
and we have the splitting forcing sequence
$$0 \lra V \lra V \oplus \AA_Y(-e_0) \lra \AA_Y(-e_0) \lra 0 \, .$$ 
Then $V \cong \PP(V') - \PP(V)$ and $\PP(V')$ is just the projective
closure of $V$.
\end{example}

\begin{proposition}
\label{solid}
Let $R$ be a normal standard graded $K$-algebra of dimension $d \geq 2$,
let $f_1, \ldots ,f_n \in R$ be primary homogeneous elements
and let $f_0$ be another homogeneous element.
Let $V$ and $V'$ be as in {\rm \ref{forcingsequence2}}.
Then the following are equivalent.
\renewcommand{\labelenumi}{(\roman{enumi})}
\begin{enumerate}
\item
$f_0 \in (f_1, \ldots ,f_n)^\soclo$.

\item
The cohomological dimension of
$\PP(V') - \PP(V)$ is $d -1 = \dim \, Y$.
\end{enumerate}
In particular, if $d=2$, then $f_0 \in (f_1, \ldots ,f_n)^\soclo$ holds
if and only if the open subscheme $\PP(V')- \PP(V)$ is not affine.
\end{proposition}
\proof
We have $\PP(V') -\PP(V) \cong D_+(R_+) \subseteq
\Proj \, R[T_1, \ldots ,T_n]/(f_1T_1 + \ldots +f_nT_n+f_0)$.
In general, every quasicoherent sheaf on an open subset
$D_+(\foa) \subseteq \Proj\, S$
is quasicoherent extendible to $\Proj\, S$ and hence of type $\tilde{M}$,
where $M$ is a graded $S$-module
(\cite[Propositions II.5.8 and II.5.15]{haralg}).
Therefore the cohomological dimensions of $D_+(\foa)$ and of $D(\foa)$
are the same.
Hence the cohomological dimensions of $\PP(V') - \PP(V)$ and of
the forcing affine-linear bundle
$D(R_+) \subseteq \Spec R[T_1, \ldots ,T_n]/(f_1T_1 + \ldots +f_nT_n+f_0)$
are the same, and the result follows from \ref{solidcd}(ii).
\qed

\begin{lemma}
\label{plus}
Let $R$ be a normal standard graded $K$-algebra of dimension $d \geq 2$,
let $f_1, \ldots ,f_n \in R$ be primary homogeneous elements
and let $f_0$ be another homogeneous element.
Let $V$ and $V'$ be as in {\rm \ref{forcingsequence2}}.
Then the following are equivalent.
\renewcommand{\labelenumi}{(\roman{enumi})}
\begin{enumerate}

\item
$f_0 \in (f_1, \ldots ,f_n)^{+{\rm gr}}$, i.e. there exists
a finite graded extension $R \subseteq R'$ such that
$f_0 \in (f_1, \ldots ,f_n)R'$.

\item
There exists a finite surjective morphism $g: Y' \ra Y$ such that
the pull back $g^*\PP(V')= \PP(V') \times_Y Y'$
has a section not meeting $g^*\PP(V)= \PP(V) \times_Y Y'$.

\item
There exists a closed subvariety $Y'' \subset \PP(V')$
not intersecting $\PP(V)$, finite and surjective over $Y$.

\item
There exists a closed subvariety $Y'' \subset \PP(V')$
not intersecting $\PP(V)$ of dimension $d-1 = \dim\, Y$.

\end{enumerate}
\end{lemma}
\proof
(i) $\Leftrightarrow $ (ii).
If $R \subseteq R'$ is finite and graded such that
$f_0 \in (f_1, \ldots ,f_n)R'$, then there exists a section
$Y' = \Proj\, R' \ra g^*\PP(V')$ which
does not meet $g^*\PP(V)$ due to \ref{trivialtwo}.
If $g:Y' \ra Y$ is such a morphism, then $g^* \O_Y(1)$ is ample
on $Y'$ and this gives the homogeneous ring $R'$.

Suppose that (ii) holds. Then the image of the section gives the closed
subvariety $Y''$ finite over $Y$. This gives (iii) and then (iv).
Suppose that (iv) holds.
The mapping $Y'' \hookrightarrow \PP(V') \ra Y$ is projective
and the fibers are zero-dimensional, since $\PP(V_y) \subset \PP(V'_y)$
meets every curve, but $Y'' \cap \PP(V) = \emptyset $.
Hence this mapping is finite and due to the assumption on the dimension
it is surjective.
So suppose that (iii) holds.
The mapping
$g : Y'=Y'' \stackrel{i}{\hookrightarrow} \PP(V') \ra Y$
is finite and surjective,
and the image of the section
$ i \times id_{Y'}  :Y'  \ra \PP(V') \times_Y Y' =g^* \PP(V') $
is disjoined to $g^* \PP(V)= \PP(V) \times_Y Y'$.
\qed

\section{Ample and basepoint free forcing divisors}
\label{ample}

Let $Z=\PP(V) \subset \PP(V')$ be the forcing divisor on $Y= \Proj\, R$
corresponding to
homogeneous forcing data $f_1, \ldots ,f_n;f_0 \in R$.
When is $Z$ ample and when is $Z$ basepoint free?
For $e_0=0$ the forcing divisor is a hyperplane section of
$\O_{\PP(V')}(1)$, and the ampleness of this invertible sheaf is by
definition the ampleness of the locally free sheaf
$\shF'(-d_0)= \pi_*\O_{\PP(V')}(1)$, see \cite[III, \S 1]{haramp} and
\cite{tightslope} for further ampleness criteria for vector bundles
and applications to tight closure problems.

Throughout this section we will assume that $K$ is algebraically closed.
The following proposition shows that the ample property is interesting
only in dimension two.

\begin{proposition}
\label{eins}
Let $K$ be an algebraically closed field and let $R$ be a normal
standard graded $K$-algebra of dimension $d$.
Let $f_1, \ldots ,f_n$ be homogeneous $R_+$-primary elements and
let $f_0$ be another homogeneous element of degrees $d_i$.
Let $V_m,V'_m$ be as in {\rm \ref{forcingsequence2}}
and let $Z=\PP(V)$ be the forcing divisor.
Then the following hold.

\renewcommand{\labelenumi}{(\roman{enumi})}
\begin{enumerate}
\item
Suppose that $f_0$ is a unit and $d_i \geq 1$ for $i=1, \ldots ,n$.
Then $Z$ is ample.

\item
If $f_0$ is not a unit, then the cohomological dimension
$cd \, (\PP(V')- Z) \geq d-2 $.

\item
If $f_0$ is not a unit and $d \geq 3$, then $Z$ is not ample.

\end{enumerate}
\end{proposition}
\proof
(i).
We may assume that $f_0=1$.
Then
$$V'_m =D_+(R_+) \subset
\Proj\, R[T_0,T_1, \ldots , T_n]/(\sum_{i=1}^nf_iT_i+T_0)
\cong \Proj\, R[T_1,\ldots,T_n] \, ,$$
where
$\deg\, T_i=e_i=-d_i + e_0$.
Thus $V'_m \cong \AA_Y(d_1-e_0) \times_Y \ldots \times_Y \AA_Y(d_n-e_0)$.
For $m=d_0$ we see that $\shF'(-d_0)$ is a sum of ample invertible sheaves,
hence $\shF'(-d_0)$ is ample due to \cite[III, Corollary 1.8]{haramp}.

(ii).
For $d=0,1$ there is nothing to show, so suppose $d \geq 2$.
The zero set $V_+ (f_0) \subset Y$ is a closed subset of
dimension $\geq d-2$. There exists a section
$V_+(f_0) \ra \PP(V')$ which does not meet $Z$.
Hence $\PP(V')-Z$ contains a projective subvariety of
dimension $d-2$, thus the inequality holds for the cohomological dimension.

(iii).
Due to (ii) the complement of $Z$ cannot be affine (it contains
projective curves), hence $Z$ is not ample.
\qed

\medskip
The forcing divisor $Z$ is basepoint free if and only if
$\O_{\PP(V')}(1)$ is generated by global sections for $e_0=0$
and this is true if and only if
$\pi_* \O_{\PP(V')}(1) = \shF'(-d_0)$ is generated by global sections.
A divisor $Z$ is called {\em semiample} (\cite[Definition 2.1.14]
{lazarsfeldpositive})
if $aZ$ is basepoint free for some $a \geq 1$.
In this case there exists a (projective) morphism
$\varphi: \PP(V') \ra \PP^N$ such that $aZ = \varphi^{-1}(H)$, where $H$
is a hyperplane section in $\PP^N$.
Then $\PP(V') -Z$ is projective over $\PP^N -H$. Schemes which are
proper over an affine scheme are called {\em semiaffine} and were
studied in \cite{goodlandman}.

\begin{lemma}
\label{semiaffine}
Let $K$ be an algebraically closed field and let
$R$ be a normal standard graded $K$-algebra.
Let $f_1, \ldots ,f_n$ be primary homogeneous elements
and let $f_0$ be another homogeneous element.
Suppose that $\PP(V') -\PP(V)$ is semiaffine.
Then $f_0 \in (f_1, \ldots ,f_n)^\soclo$ if and only if
$f_0 \in (f_1, \ldots, f_n)^{+{\rm gr}}$.
\end{lemma}
\proof
Due to \cite[Corollary 5.8]{goodlandman} the cohomological dimension of
a semiaffine scheme equals the maximal dimension of a closed
proper subvariety.
Thus $f_0 \in (f_1, \ldots , f_n)^\soclo$ implies via \ref{solid} that
there exists a projective subvariety $Y' \subset \PP(V')$
of dimension $ \dim \, Y$ which does not meet $\PP(V)$.
Therefore $f_0 \in (f_1 , \ldots ,f_n)^{+{\rm gr}}$ due to \ref{plus}.
\qed

\medskip
The condition in the following corollary
is usefull only for $\dim \, R=2$.
We will apply this in section \ref{higherrank}.

\begin{corollary}
\label{pullbackample}
Let $K$ denote an algebraically closed field and let
$R$ be a normal standard graded $K$-algebra.
Let $f_1, \ldots ,f_n$ be primary homogeneous elements
and let
$Z= \PP(f_1, \ldots ,f_n) \subset \PP(f_1, \ldots, f_n;f_0)$
be the corresponding bundles on $Y=\Proj \, R$.
Suppose that the pull back $Z|_Z$ is ample.
Then
$f_0 \in (f_1, \ldots ,f_n)^\pasoclo$ if and only if
$f_0 \in (f_1, \ldots ,f_n)^{+{\rm gr}}$.
\end{corollary}
\proof
The theorem of Zariski-Fujita (see \cite[Remark 2.1.18]{lazarsfeldpositive})
asserts that
$Z$ is semiample. Then the complement of $Z$ is semiaffine
and the result follows from \ref{semiaffine}.
\qed

\begin{corollary}
\label{globalsection}
Let $R$ be a normal standard graded $K$-algebra and
let $f_1, \ldots ,f_n$ be primary homogeneous elements
and let $f_0$ be another homogeneous element.
Suppose that the locally free sheaf $\shF'(-d_0)$
is generated by global sections.
Then $f_0 \in (f_1, \ldots ,f_n)^\soclo$ if and only if
$f_0 \in (f_1, \ldots, f_n)^{+{\rm gr}}$.
\end{corollary}
\proof
Since $\shF'(-d_0)$ is generated by global sections we know that
the forcing divisor is basepoint free, hence
$\PP(V') - \PP(V)$ is semiaffine.
Hence the result follows from \ref{semiaffine}.
\qed

\medskip
Note that the results \ref{semiaffine} - \ref{globalsection}
yield in characteristic zero 
in fact the stronger result that
$f_0 \in (f_1, \ldots ,f_n)^\soclo$ holds if and only if already
$f_0 \in (f_1, \ldots , f_n)$ holds.
The following corollary was proved for tight closure
in \cite[Theorem 2.2]{smithgraded}
using differential operators in positive characteristic.
Our version proves the same result for solid closure.

\begin{corollary}
\label{degreevanish}
Let $K$ denote an algebraically closed field and
let $R$ be a normal standard graded $K$-algebra and
let $f_1, \ldots ,f_n$ be primary homogeneous elements of degrees
$d_i$.
Let $f_0$ be another homogeneous element of degree
$d_0 \leq \min_{i} d_i $.
Then $f_0 \in (f_1, \ldots ,f_n)^\soclo$ is only possible if
$f_0 \in (f_1, \ldots, f_n)$.
\end{corollary}
\proof
Set $e_0=0$.
Then $e_i=d_0 -d_i \leq 0$ and we have a surjection
$\O_Y(-e_1) \oplus \ldots \oplus \O_Y(-e_n) \oplus \O_Y \ra \shF'(-d_0) \ra 0$.
Since the $\O_Y(a)$ for $a \geq 0$ are generated by global sections,
we have also a surjection $\O_Y^k \ra \shF'(-d_0) \ra 0$.
Therefore $\shF'(-d_0)$ is generated by global sections and we have a closed
embedding $V' \hookrightarrow Y \times \AA^k$.

Suppose that $f_0 \in (f_1, \ldots ,f_n)^\soclo$. Then by \ref{globalsection}
we know that there exists a subvariety $Y' \subset \PP(V')$
of dimension $\dim \, Y$ not meeting the forcing divisor $Z$.
We may consider $Y' \subset V_+(T_0-1) \subset V'$, since
$V_+(T_0-1)$ is isomorphic to $\PP(V')-\PP(V)$ via the cone mapping
(see the proof of \ref{forcingsequence2}(ii)).
All together we get a closed embedding
$Y' \hookrightarrow Y \times \AA^k$. Since $Y'$ is a projective variety,
this factors through $Y \times \{P\}$, where $P \in \AA^k$ is a closed point,
and so $Y' \cong Y \times \{P\} \cong Y$, since $K$ is algebraically closed.
Hence we get a section. 
\qed

\medskip
Even if the forcing divisor is not basepoint free, the existence of
linearly equivalent effective divisors has consequences on
the existence of closed subvarieties and hence on the existence of
finite solutions (in the sense of \ref{plus} (iii) or (iv))
for the tight closure problem.
See also Proposition \ref{nupositiv}.

\begin{proposition}
\label{effectivevertreter}
Let $R$ be a normal standard graded $K$-algebra such that
$Y=\Proj\, R$ is a smooth variety.
Let $f_1, \ldots ,f_n$ be homogeneous primary elements
and let $f_0$ be another homogeneous element.
Suppose that there exists a positive {\rm(}effective $\neq 0${\rm )}
divisor $ L \subset Y$ such that for some $a \geq 1$
the divisor $a \PP(V)- \pi^*L $ is linearly equivalent
to an effective divisor.
Then there exists a linearly equivalent effective divisor
$D \sim a\PP(V)$ with the property that
the cohomological dimension of $\PP(V') -{\rm supp}\,  D$ is smaller
than the {\rm(}cohomological{\rm)} dimension of $Y$.

If $Y' \subseteq \PP(V')$ is finite
and surjective over $Y$ and disjoined to $\PP(V)$
{\rm(}as in {\rm \ref{plus}(iii))},
then $Y'$ must lie on the support of $D$.
\end{proposition}
\proof
Let $a\PP(V)- \pi^*L \sim D' $ be effective, hence $a\PP(V) \sim D=D' +\pi^*L$.
The divisor $D'$ cuts out a hyperplane on every fiber,
hence it is also a projective subbundle.
Since a projective bundle minus a dominant effective divisor is relatively
affine over the base we see that
$\PP(V') - {\rm supp}\,D$ is affine over $Y-{\rm supp}\,L$.
But the cohomological dimension of $Y-{\rm supp}\,L$ is smaller than the
dimension of $Y$ due to the theorem of Lichtenbaum
(\cite[Corollary 3.2]{hartshornecohdim}),
hence this is also true for $\PP(V') -{\rm supp}\,D$.

Now suppose that $Y'$ is finite and surjective over $Y$
and $Y' \cap \PP(V) = \emptyset$.
Then we have from intersection theory the identities 
$0=a Y' . \PP(V)= Y'. (D'+\pi^*L) =Y'.D' + Y'.\pi^*L$. The second
summand is a positive cycle, since $Y'$ dominates $Y$.
Hence $Y'.D'$ cannot be effective and the intersection
of $Y'$ and $D'$ must be improper,
so $Y' \subset {\rm supp}\, D' \subset {\rm supp}\, D$.
\qed

\bigskip
For the rest of this paper we will restrict to the situation
where $K$ is an algebraically closed field
and $R$ is a two-dimensional normal standard graded
$K$-algebra. Then $Y= \Proj \, R$ is a smooth projective curve
with hyperplane section $H$.
Homogeneous primary elements $f_1, \ldots ,f_n,f_0$ yield the
projective bundle $\PP(V')=\PP(f_1, \ldots ,f_n;f_0)$ of rank $n-1$
over the curve (and of dimension n) together with the
forcing divisor $Z=\PP(f_1, \ldots ,f_n)$.
Now $f_0 \in (f_1, \ldots , f_n)^\soclo$ holds if and only if the complement
of the forcing divisor is not affine.

If the complement of the forcing divisor is affine
(i.e. $f_0 \not\in (f_1, \ldots ,f_n)^\soclo$),
then $\Gamma(\PP(V')- \PP(V), \O_{\PP(V')})$ is a finitely generated
$K$-algebra of dimension $n$.
It follows that some multiple of the forcing divisor $\PP(V)$
defines a rational mapping to some projective space
such that the dimension of the image is
$n$. This means by definition that
$\PP(V)$ is {\em big} (has maximal Iitaka-dimension,
see \cite[Definition 2.2.1]{lazarsfeldpositive}).

It is sometimes possible to establish
$f_0 \in (f_1, \ldots ,f_n)^\soclo$ by showing that the forcing divisor is not big.
The following proposition deals with equivalent conditions for bigness
in our situation.

\begin{proposition}
\label{bigaffine}
Let $R$ be a normal two-dimensional standard graded $K$-algebra.
Let $f_1, \ldots ,f_n$ be homogeneous primary elements
and let $f_0$ be another homogeneous element.
Let $Z=\PP(V) \subset \PP(V')$ denote the forcing divisor.
Then the following are equivalent.

\renewcommand{\labelenumi}{(\roman{enumi})}
\begin{enumerate}
\item
There exists a positive divisor $L \subset Y$ such that
for some $a \geq 1$ the divisor
$a Z - \pi^*L$ is equivalent
to an effective divisor.

\item
There exists a linearly equivalent effective divisor
$D \sim a Z$ {\rm(}$a \geq 1${\rm )} such that
$\PP(V')- {\rm supp} \, D$ is affine.

\item
The forcing divisor $Z$ is big.

\end{enumerate}

\end{proposition}
\proof
(i) $\Rightarrow$ (ii) follows from \ref{effectivevertreter}.
Suppose that (ii) holds, let $X= \PP(V')$ and let 
$s \in \Gamma(X,\O_X(aZ))$ be a section such that $X_s$ is affine.
The topology of $X_s$ is generated by
subsets $X_t \subseteq X, \, t \in \Gamma(X,\O_X(bZ))$, $ b \geq 1$,
see \cite[Th\'{e}or\`{e}me 4.5.2]{EGAII}.
Therefore the rational mapping defined by $aZ$ is an isomorphism
on $X_s$ and the image has maximal dimension, hence $Z$ is big
(and (ii) $\Rightarrow $ (iii)).
On the other hand, 
if $\emptyset \neq V \subset X_s$ is an affine subset which
does not meet the fiber over a point $P \in Y$,
then there exists also $t \in \Gamma(X,\O_X(bZ))$
such that $ \emptyset \neq X_t \subseteq V$.
Therefore
$bZ+(t) = \sum_i a_i D_i$, $a_i >0$ is an effective divisor
and $\PP(V'_P)$ is one of the $D_i$ (hence (ii) $\Rightarrow $ (i)).
(iii) $\Rightarrow $ (ii).
If $Z$ is big, then for some $a \geq 1$
the multiple $aZ$ defines a mapping which is birational with its image.
Therefore the mapping induces an isomorphism on an open affine subset
$X_s \cong D_+(s)$, $s \in \Gamma(X,\O_X(aZ))$.
\qed

\begin{remark}
\label{bignef}
A numerically effective divisor $Z$ is big if and only
if its top self intersection number $Z^n$ is $>0$, see
\cite[Theorem VI.2.15]{kollar} or \cite[Theorem 2.2.14]{lazarsfeldpositive}.
The top self intersection of the forcing divisor $\PP(f_1, \ldots ,f_n)
\subset \PP(f_1, \ldots ,f_n;f_0)$ corresponding to
forcing data in a two dimensional normal graded domain $R$
is $(d_1 + \ldots +d_n- (n-1)d_0) \deg\, H$,
where $H$ is the hyperplane section on $Y= \Proj \, R$.
This follows from \ref{buendelalsproj}(iii).
\end{remark}

\section{Ruled surfaces and forcing sections}
\label{ruled}

In this section $K$ denotes an algebraically closed field
and $R$ denotes a standard graded two-dimensional $K$-algebra
and we consider the tight closure of homogeneous
parameters $f_1$ and $f_2$.
The construction of projective bundles and subbundles
from homogeneous elements described in section \ref{graded}
leads in this setting to ruled surfaces over
the curve $\Proj \, R$ together with a forcing section.
It is known that
the so-called vanishing theorem
$(f_1,f_2)^\soclo = (f_1,f_2) +R_{\geq \deg\, f_1 + \deg \, f_2}$
holds for $\Char (K)=p =0$ or $p >>0$
(see \cite[Theorem 4.3]{hunekesmithkodaira}), 
and we will prove this result in our geometric setting.

\begin{corollary}
\label{forcingrule}
Let $K$ denote an algebraically closed field
and let $R$ be a standard graded two-dimensional normal $K$-algebra,
let $f_1,f_2$ be homogeneous parameters of degrees $d_1$ and $d_2$ and
let $f_0$ be another element of degree $d_0$.
Set $\leadno = d_1 +d_2-d_0 $. Let $e_1,e_2,e_0$ be integers such that
$e_i+d_i =m $ and let $V_m$ and $V'_m$ as in \ref{forcingsequence2}
and set $Y= \Proj \, R$. Then the following hold. 

\renewcommand{\labelenumi}{(\roman{enumi})}
\begin{enumerate}

\item
$\PP(V')$ is a ruled surface and $\PP(V) \subset \PP(V')$
is a section, called the forcing section.

\item 
We have
$\Proj \,R[T_1,T_2]/(f_1T_1+f_2T_2) \supset D_+(R_+) =V_m
\cong \AA_Y( \leadno -e_0)$.
In particular, the exact forcing sequence is
$$ 0 \longrightarrow \AA_Y( \leadno -e_0) \longrightarrow V'_m
\longrightarrow \AA_Y(-e_0) \longrightarrow 0 \, .$$

\item
We have ${\rm Det} \, V' \cong \AA_Y( \leadno -2e_0)$.

\item
The normal bundle for the embedding
$Y\cong \PP(V) \subset \PP(V')$
is $\AA_Y( -\leadno  )$.

\item
The self intersection number of the forcing section
$Y \cong \PP(V) \hookrightarrow \PP(V')$
is $ \leadno \deg\, H$, where
$H$ is the hyperplane section corresponding to $\O_Y(1)$.
\end{enumerate}
\end{corollary}

\proof
(i)
follows from \ref{forcingsequence2}.
(ii)
The homomorphism
$R[T_1,T_2]/(f_1T_1+f_2T_2) \ra R[T]$ given by
$T_1 \mapsto f_2T,\, T_2 \mapsto -f_1T$ is homogeneous for
$\deg\, T=e_1-d_2=e_2-d_1$ and induces an isomorphism
on $D_+(R_+)$.
Since $\deg \, T = e_1-d_2= m-d_1-d_2= e_0 -\leadno  $
the corresponding line bundle is $\AA_Y(\leadno-e_0)$. (iii) follows.

The normal bundle for the embedding on $\PP(V) \cong Y$ is
$N = \AA_{\PP(V)}(-1) \otimes \AA_Y(-e_0)$ due to \ref{forcingsequence2}(iv).
Furthermore,
$V_m=  \AA_Y( \leadno -e_0)$ on $Y$ and
$V_m=\AA_{\PP(V)}(+1)$ is the tautological line bundle
for $\PP(V) \cong Y$. This yields
$N = \AA_Y(-\leadno +e_0) \otimes \AA_Y(-e_0)=\AA_Y( -\leadno  )$.
Its sheaf of sections is
$\O_Y( \leadno )$ and its degree is the self intersection number,
hence (iv) follows.
\qed

\begin{remark}
\label{forcingsheaf3}
The corresponding sequence of sheaves are
$$ 0 \longrightarrow \O_Y(-e_0) \longrightarrow \shF'(-m)
\longrightarrow \O_Y( \leadno -e_0) \longrightarrow 0 \, $$
for the linear forms $\shF'(-m)$ and
$$ 0 \longrightarrow \O_Y(- \leadno +e_0) \longrightarrow \shR'(m)
\longrightarrow \O_Y(e_0) \longrightarrow 0 \, $$
for the relations $\shR'(m)$ .
These extensions are classified by
$H^1(Y, \O_Y(- \leadno))$ for $e_0=0$, where the elements
$f_1,f_2;f_0$ correspond to the cohomology class
$f_0/f_1f_2$.
\end{remark}

The following proposition gives a criterion for tight closure
in the two-dimensional parameter case in terms of ampleness of the
forcing divisor.

\begin{theorem}
\label{ampelcrit}
Let $R$ be a two-dimensional standard graded normal $K$-algebra and let
$Y= \Proj\, R$.
Let $f_1,f_2$ be homogeneous parameters and let
$f_0$ be another homogeneous element.
Let $s:Y \ra \PP(V')$ be the corresponding forcing section,
$Z=\PP(V)=s(Y)$.
Then the following are equivalent.

\renewcommand{\labelenumi}{(\roman{enumi})}
\begin{enumerate}

\item 
$f_0 \not\in (f_1,f_2)^\pasoclo$.

\item
$\PP(V')-Z$ is affine.

\item
The forcing divisor $Z$ on $\PP(V')$ is ample.
\end{enumerate}
\end{theorem}

\proof
We know the equivalence (i)$ \Leftrightarrow$ (ii) from proposition \ref{solid}
so we have to show the equivalence (ii)$ \Leftrightarrow$ (iii).
If $Z$ is ample, then its complement is affine.
If $\PP(V')-Z$ is affine, then it does not
contain projective curves.
Furthermore, there exist global functions on $\PP(V')-Z$ which are
not constant.
Thus $aZ$, $a \geq 1$,
is linearly equivalent with an effective divisor not containing $Z$.
Hence the self intersection number is positive and the
criterion of Nakai yields that $Z$ is ample.
\qed

\begin{corollary}
\label{degreecrit1}
Let $K$ be an algebraically closed field and let $R$ be a normal
two-dimensional
standard graded $K$-algebra and let $f_1,f_2$
be homogeneous parameters of degrees
$d_1,d_2$.
Then $R_{ \geq d_1 + d_2} \subseteq (f_1,f_2)^\pasoclo$.

If the characteristic of $K$ is zero,
then
$(f_1,f_2)^\pasoclo= (f_1,f_2)+ R _{\geq d_1 + d_2}$
\end{corollary}
\proof
Let $\deg\, f_0 \geq d_1+d_2$.
Then $\leadno \leq 0$ and the self intersection of
$Z=\PP(V) \subset \PP(V')$
is not positive. Hence $Z$ is not ampel and $f_0 \in (f_1,f_2)^\pasoclo$
due to theorem \ref{ampelcrit}.

Now let $f_0 \in (f_1,f_2)^\pasoclo$, but $f_0 \not\in R_{\geq d_1+d_2}$.
Then the self intersection is positive, but
the forcing divisor $Z$ is not ample.
Thus there must exist a curve $C \subset \PP(V')$ disjoint to $Z$.
By \ref{plus} it follows that $f_0 \in (f_1,f_2)^+$,
so \ref{finiteextension}(ii) gives the result in characteristic $0$.
\qed

\medskip
To prove the result of \ref{degreecrit1}
also in positive characteristic, we need the notion
of a normalized section and of
the so-called $e$-invariant of a ruled surface.
From the point of view of the forcing divisor it is technically
convenient to introduce the normalizing number.

\begin{definition}
\label{normalnumberdef}
Let $W$ be a vector bundle over a smooth projective curve $Y$ and let
$\pi: \PP(W) \ra Y$ be the projective bundle. Let
$D$ be a divisor on $\PP(W)$.
We say that $D$ is {\em normalized} if $D$ has an effective representative,
but for every divisor $\fod $ of $Y$
of negative degree the divisor $D + \pi^* \fod $ does not
have an effective representative.

For any divisor $Z$ on $\PP(W)$ we call the number
$\nonu $ characterized by the fact that there exists
a divisor $\fod $ on $ Y$ of degree  $- \nonu$
such that $Z + \pi^* \fod$ is normalized the
{\em normalizing number} of $Z$.

If $Z=\PP(V) \subset \PP(V')$
is the forcing divisor of a tight closure problem
in a two-dimenional normal standard graded $K$-algebra,
then we call the normalizing number of $Z$
also the normalizing number of the problem or of the forcing data.
\end{definition}

\begin{remarks}
\label{einvariantedef}
Recall that a locally free sheaf $\shE$ on a smooth projective curve $Y$
is called normalized if
$H^0(Y,\shE) \neq 0$,
but $H^0(Y, \shE \otimes \shL)=0$ for every invertible sheaf $\shL$
of negative degree.
If $\shE$ is normalized, then a section $0 \neq s \in H^0(Y , \shE) $
yields a mapping $s^\dual : \shE^\dual \ra \O_Y$, which must be surjective.
Hence the kernel is locally free and we get
a short exact sequence
$0 \ra \O_Y \stackrel{s}{\ra} \shE \ra \shF \ra 0$ where
$\shF$ is locally free
and where the corresponding projective subbundle $\PP(\shF) \subset \PP(\shE)$
is normalized.

On the other hand, if
$\PP(\shF) \subset \PP(\shE)$ is normalized, then
the sheaf
$\shE \otimes \ker (\shE \ra \shF )^\dual$ is normalized.

Suppose that $\PP(W)=\PP(\shE)$ is a ruled surface.
Then a normalized subbundle of codimension one is the same as
a {\em normalized section}. Recall that the
{\em $e$-invariant} of a ruled surface is defined
by $e=- C_0^2 $, where $C_0$ is a normalized section,
and that $e = - \deg \, \shE$ for normalized $\shE$.

If $\shF'(-m)$ is the sheaf of linear forms coming from a tight
closure problem, then $\nonu$
is also
characterized by the property that there exists
an invertible sheaf $\shL$ of degree $\, e_0 \deg\, H- \nonu \,$
such that $\, \shF'(-m) \otimes \shL\,$ is normalized.
If $H^0(Y, \shF'(-m) \otimes \shM)=0$ for every invertible
sheaf $\shM$ of negative degree,
then $\nonu \leq (m-d_0) \deg \, H$.
\end{remarks}

\begin{lemma}
\label{lunde}
Let $f_1,f_2,f_0 \in R$ as in {\rm \ref{forcingrule}}.
Let $e$ denote the $e$-invariant of the ruled surface
$\PP(f_1,f_2;f_0)$,
let $\leadno=d_1+d_2-d_0$ and let $\nonu$ be the normalizing number of
$Z=\PP(V)$. Then
$$e=2 \nonu - \leadno \deg\, H\, \, \mbox{ and }\,
\,\nonu = \frac{\leadno \deg\, H+e}{2} \, .$$
\end{lemma}

\proof
Let $\fod $ be a divisor on $Y$ of degree $-\nonu$ such that
$Z + \pi^* \fod $ is normalized and let
$C_0 \sim Z+ \pi^* \fod$ be effective.
Hence $C_0$ is a normalized section.
Numerically we have $C_0 = Z- \nonu F$ ($F =\,$ fiber)
and therefore
$-e=C_0^2= Z^2-2 \nonu =  \leadno  \deg\, H -2 \nonu  $.
\qed

\medskip
Knowing the $e$-invariant of a ruled surface one may
characterize the divisors which are ample.
We recall this only for sections.
\begin{lemma}
\label{amplerulecrit}
Let $S$ be a ruled surface with $e$-invariant $e$
and let $D$ be a section.
Then the following hold.

\renewcommand{\labelenumi}{(\roman{enumi})}
\begin{enumerate}

\item
Suppose that $e \geq 0$. Then $D$ is ample if and only if $D^2 >e$.

\item
Suppose that $e <0$ and that $\Char \, K =0 $ or $p >> 0$.
Then $D$ is ample.

\end{enumerate}
\end{lemma}

\proof
Let $C_0$ be a normalized section and write $D=C_0+bF$ where
$ b \geq 0$.
Then $D^2=C_0^2+2b$, thus $b=1/2(D^2-C_0^2)$.
If $e \geq 0$, then the criterion \cite[Proposition V.2.20]{haralg}
says $b > e$, thus $1/2(D^2-C_0^2) > e=-C_0^2$ and this gives
$D^2 > e$.

Let $e<0$. If the characteristic is zero, then \cite[Proposition 2.21]{haralg}
gives the result. If the characteristic is positive,
then \cite[Excercise 2.14]{haralg}
yields the condition
$b> (e/2 + (g-1)/p)$, and this is true for $p>> 0$.
\qed

\medskip
To apply this criterion on ruled surfaces arising from
forcing equation, we need to know something about
the $e$-invariant of them.

\begin{lemma}
\label{einvariante1}
Let $f_1,f_2,f_0 \in R$ be as in {\rm \ref{forcingrule}}.
Let $\shE=\shF'(-m)$ be the sheaf of linear forms of $V'_m$.
Then the following hold.

\renewcommand{\labelenumi}{(\roman{enumi})}
\begin{enumerate}

\item
If
$e_0 \leq 0$, then $H^0(Y,\shE) \neq 0$.

\item
Let $e_0 \geq 0$.
Let $\shL$ be an invertible sheaf on $Y$ of negative degree $-k$.
Then $H^0(Y,\shE \otimes \shL)=0$
for $k > ( \leadno -e_0) \deg H $.
If $f_0 \not\in (f_1,f_2)$, then this is also true for
$k \geq ( \leadno -e_0) \deg H $.

\item
Let $ \leadno  \leq 0$ and $e_0=0$.
Then $\shE$ has global sections $\neq 0$, but
$\shE \otimes \shL$ does not have for $\deg \shL <0$
{\rm (}i.e. $\shE = \shF'(-d_0)$ is normalized{\rm )}.

\end{enumerate}
\end{lemma}
\proof

(i).
If $e_0 \leq 0$, then
$0 \neq H^0(Y,\O_Y(-e_0)) \subseteq H^0(Y,\shE)$ by the forcing sequence
and
$\shE$ has global sections $\neq 0$.

(ii)
We consider the exact forcing sequence of the sheaf of linear forms
$$0 \longrightarrow \O_Y(-e_0) \longrightarrow \shE
 \longrightarrow \O_Y( \leadno -e_0) \longrightarrow 0\, .$$
We tensorize with $\shL$ and get
the cohomology sequence
$$0 \longrightarrow H^0(Y,\O_Y(-e_0) \otimes \shL)
\longrightarrow H^0(Y,\shE \otimes \shL)
\longrightarrow H^0(Y,\O_Y( \leadno -e_0) \otimes \shL) \, .$$
Because of
$e_0 \geq 0$ and $\deg \, \shL < 0$ we have on the left hand side
and because of
$k > ( \leadno -e_0) \deg \,H $ we have on the right hand side
an invertible sheaf of negative degree, thus
the result follows.

If $k=( \leadno -e_0) \deg \,H $, then on the right hand side
we have an invertible sheaf of degree zero.
If it is not trivial, then it
has no global sections.
Otherwise the sheaf is $\O_Y$ and the cohomology sequence is
$$0 \longrightarrow H^0(Y,\shE \otimes \shL) \longrightarrow
H^0(Y,\O_Y)=K
\longrightarrow H^1(Y,\O_Y(-e_0) \otimes \shL) \, .$$
Suppose that $H^0(Y,\shE \otimes \shL) \neq 0$. Then
this maps surjective on $K$ and $K$ maps to zero.
But then the short exact sequence splits, which
contradicts the assumption $f_0 \not \in (f_1,f_2)$.

(iii) follows from (i) and (ii).
\qed

\begin{corollary}
\label{einvariante2}
Let $f_1,f_2,f_0$ be as in {\rm \ref{forcingrule}}
and let $e$ denote the $e$-invariant of
$\PP(f_1,f_2 ;f_0)$
and let $\nonu$ denote the normalizing number.
Then the following hold.

\renewcommand{\labelenumi}{(\roman{enumi})}
\begin{enumerate}

\item
We have $\nonu \geq 0$ and $e \geq -\leadno \deg \, H$.

\item
Let $ \leadno  >0$.
Then $\nonu \leq \leadno \deg \, H$ and
$\nonu < \leadno \deg \, H$ if $f_0 \not\in (f_1,f_2)$
{\rm(}and
$e \leq \leadno \deg\, H$ and $e < \leadno \deg\, H$ if
$f_0 \not\in (f_1,f_2)${\rm )}.

\item
Let $ \leadno  \leq 0$.
Then $\nonu=0$ and $e= -\leadno \deg \, H \geq 0$.
Thus $\nonu >0$ {\rm(}or $e < 0${\rm )} implies that $ \leadno  >0$.
\end{enumerate}
\end{corollary}

\proof
The statements on $e$ follow from the statements on $\nonu$ by \ref{lunde}.
Fix $e_0=0$ and let $\shE = \shF'(-d_0)$
be the corresponding sheaf of linear forms.
In determining $\nonu$
we have to look for which invertible sheaves $\shL$ there exist sections
$0 \neq s \in \Gamma(Y, \shE \otimes \shL)$.
So the results follow all from the corresponding statements in
\ref{einvariante1}.
(ii). From \ref{einvariante1}(ii) we get that
$H^0(Y, \shF'(-d_0) \otimes \shL)=0$ for every
invertible sheaf $\shL$ of degree $-k$, $k > \leadno \deg \, H$
(so $\shL$ is automatically of negative degree), hence
$\nonu \leq \leadno \deg \, H$.
\qed

\begin{corollary}
\label{degreecrit2}
Let $K$ be an algebraically closed field of characteristic $p$
and let $R$ be a two-dimensional normal
standard graded $K$-algebra and let
$f_1,f_2$ be homogeneous parameters of degrees
$d_1,d_2$.
Then for $p=0$ and for $p >> 0$ we have
$(f_1,f_2)^\pasoclo= (f_1,f_2)+ R _{\geq d_1+d_2}$.
\end{corollary}

\proof
We have proved the inclusion $\supseteq $ in \ref{degreecrit1}.
Thus suppose that
$f_0 \not\in (f_1,f_2)$ and $ f_0 \not\in  R _{\geq d_1+d_2}$.
Due to \ref{ampelcrit} we have to show that the forcing divisor
$Z=\PP(f_1,f_2) \subset \PP(f_1,f_2;f_0)$
is ample.
The self intersection number of $Z$
is positive due to the second assumption.
If $e < 0$,
then the statement follows from \ref{amplerulecrit}(ii)
for $p>>0$.
If $e \geq 0$, then \ref{einvariante2}(ii) shows that
$\leadno \deg\, H > e$ and again \ref{amplerulecrit}(i) gives ampleness.
\qed

\section{Examples}
\label{examples}
We give some examples of ruled surfaces and their forcing
sections arising from tight closure problems.
Let $R$ be a standard graded normal two-dimensional $K$-algebra,
$Y= \Proj \, R$,
and let $f,g$ be homogeneous parameters and $h$ homogeneous
of degrees $d_1,d_2,d_0$. Set $\leadno = d_1+d_2-d_0$.

\begin{example}
\label{spaltet}
Suppose that $h \in (f,g)$.
Then the forcing sequence splits and
$$V'_m \cong \AA_Y( \leadno -e_0) \times \AA_Y(-e_0)$$
This is normalized for $e_0= \max (0,  \leadno  )$
and then isomorphic to
$\AA_Y \times \AA_Y(-| \leadno |)$ and
$\PP(V')$ is the projective closure of the line bundle
$\AA_Y(-|\leadno|)$. The forcing section is either the zero section
of the line bundle or the closure section.
The $e$-invariant is $|\leadno|\deg \, H$.
\end{example}

\begin{example}
\label{zerlegbar}
Suppose that $h=1$. Then
$$\Proj\, R[T_0,T_1,T_2]/(fT_1+gT_2+T_0) \cong \Proj\, R[T_1,T_2] \, .$$
and
$V'_m= \AA_Y(-e_1) \times \AA_Y(-e_2)$ is decomposable.
The forcing sequence is
$$0 \longrightarrow \AA_Y(d_1+d_2-e_0) \longrightarrow
\AA_Y(d_1-e_0)\times \AA_Y(d_2-e_0) \longrightarrow \AA_Y(-e_0)
\longrightarrow 0 \, ,$$
where $v \mapsto (gv,-fv)$ and $(a,b) \mapsto -(af+bg)$.
The sequence is normalized for $e_0= \max(d_1,d_2)$.
Let $d_1 \geq d_2$. Then the normalized sequence is
$$0 \longrightarrow \AA_Y(d_2) \longrightarrow
\AA_Y \times \AA_Y(d_2-d_1) \longrightarrow \AA_Y(-d_1)
\longrightarrow 0 \, .$$
The $e$-Invariant of the ruled surface is $e=|d_2-d_1|\deg\, H$.
\end{example}

\begin{example}
\label{graph}
Let $f$ and $g$ be of the same degree $d$ and let $h=1$.
This yields the trivial ruled surface $Y \times \PP^1$ and the 
forcing section is the graph of the  meromorphic function $-f/g$.
For let  $e_0=d$. Then we have the forcing sequence
$$0 \longrightarrow \AA_Y(d) \stackrel{g,-f} {\longrightarrow}
\AA_Y \times \AA_Y 
\stackrel{fT_1+gT_2}{\longrightarrow} \AA_Y(-d)
\longrightarrow 0 \, ,$$
and $fT_1+gT_2=0$ is equivalent with $T_1/T_2= -g/f$.
If $d \geq 1$, then $h \not\in (f,g)^\pasoclo$, and so
we see via tight closure that the complement of the graph
of a non constant meromorphic function is affine.
\end{example}

\begin{example}
\label{projektivegerade}
Let $R=K[x,y]$, thus $Y=\PP_K^1$.
The ruled surface $\PP(f,g;h)$ must be a Hirzebruch surface
$\PP(\AA_{\PP_K^1} \times \AA_{\PP^1_K} (k))$
and we have to determine the number $k \leq 0$.

Let $(a_1,a_2,a_3)$ and $(b_1,b_2,b_3)$ be a basis of
homogeneous relations for all the relations
for $(f,g,h)$.
Let $e_i$ and $e_i'$ be its degrees and suppose that $e_i \geq e_i'$.
Set $ k= e_i'- e_i$.
The homomorphism
$$K[x,y,T_1,T_2,T_3]/(fT_1+gT_2+hT_3) \stackrel{\psi }{\lra} K[x,y,U,S] $$
given by $T_i \mapsto a_iU+b_iS,\, i=1,2,3$, is well defined
and is homogeneous for
$\deg \, T_i=e_i,\, \deg\, U=0$, $ \deg\, S=e_i -e_i'=-k $.
It induces a mapping
$$ \Proj\, K[x,y,U,S] \supseteq D_+(x,y)
=\AA_{\PP^1_K} \times \AA_{\PP^1_K}(k)  \lra V'_m \, .$$
We claim that this is an isomorphism.
Since $(-g,f,0)$ is a relation, there exist $r,s \in R$
such that
$(-g,f,0)= r(a_1,a_2,a_3) + s(b_1,b_2,b_3)$.
Since $(0,-h,g)$ and $(-h,0,f)$ are also relations it follows that
$(g,f) \subseteq (a_3,b_3)$. Hence $(a_3,b_3)$ is $R_+$-primary
and $a_3$ and $b_3$ do not have a common divisor.
Therefore there exists $t \in R$ such that
$r=tb_3,\, s=-ta_3$, hence we may write $f=t(b_3a_2-a_3b_2)$.
The mapping $\psi$ is locally on $D(f)$ given by a linear transformation
$R_f[T_2,T_3] \ra R_f[U,S]$ and its determinant is
$b_3a_2-a_3b_2$, which is a unit in $R_f$.
The same is true on $D(g)$, so the induced mapping is an isomorphism.

The forcing sequence is
$$0 \lra \AA_{\PP_K^1} (\leadno -e_3) \lra 
\AA_{\PP^1_K} \times \AA_{\PP^1_K}(e_3'-e_3)
\stackrel{a_3U+b_3S}{\lra} \AA_{\PP_K^1}(-e_3) \lra 0
\, .$$
It follows that $\leadno =e_3 + e_3'$.

Let $f=x^2,g=y^2$. For $h=xy$ we have $e_3=e_3'=1$, but for
$h=x^2$ we have $e_3=2, \, e_3'=0$.
\end{example}

Let $f,g,h \in R_d$ be homogeneous of the same degree $d$ and
suppose that $f,g$ are parameters and that $h \not\in (f,g)$.
Let $E \subseteq \Gamma(Y, \O_Y(d))$ denote the linear system
spanned by $f,g,h$ and set $e_i=0$.
Then the sequence from \ref{buendelalsproj}(ii) is
$$0 \lra V'_d \lra \AA_Y^3 \lra \AA_Y(-d) \lra 0 \,  .$$
The mapping sends
$(P,t_1,t_2,t_3) \mapsto t_1f(P) +t_2g(P) +t_3h(P)$ and this is
zero if and only if the point $P$ lies on the divisor
defined by the section $t_1f +t_2g+t_3h$.
Therefore
$$ \PP(V')= \{ (P,D): \,\,  P \in D, \, \, D \in E \}  \, $$
and the ruled surface $\PP(V')$ is the incidence variety
to the linear system $E$.

Now let $Y =V_+(F) \subset \PP_K^2 = \Proj\, K[x,y,z]$
be a smooth curve and consider the linear system of lines.
Thus the ruled surface associated to the vector bundle
$$V'_1=D_+(x,y,z) \subset \Proj\, K[x,y,z]/(F)[T_1,T_2,T_3]/(xT_1+yT_2+zT_3) \,$$
($\deg \, T_i=0$)
consists of the pairs $(P,L)$, where $L$ is a line through $P \in Y$.
Suppose that $x,y$ are parameters for the curve. Then the point $Q=(0,0,1)$
does not belong to $Y$.
The forcing section maps a point $P \in Y$ to the line passing through
$P$ and $Q$, since this is the line given by $T_3=0$.
The self intersection number of the forcing section is
$\deg \,F$. If $\deg \, F \geq 2$, then $z \not\in (x,y)^\pasoclo$ and
the forcing section is ample.

\begin{example}
Let $F \in K[x,y,z]$ be a homogeneous polynomial of degree three
such that $Y= \Proj \, K[x,y,z]/(F)$ is an elliptic curve and
suppose that $x,y$ are parameters.
Consider $\PP(x,y;z)$ and set $e_1=e_2=e_3=0$, and let $\shE$
be the sheaf of linear forms for this grading.
The global linear forms $T_1,T_2,T_3$ are a basis for $H^0(Y, \shE) $.
A global linear form $s=a_1T_1+a_2T_2+a_3T_3$ ($a_i \in K$)
belongs to $\Gamma(Y, \shE \otimes \O_Y(-P))$ if and only if
$s|V'_P =0$, and this is the case if and only if
$(a_1,a_2,a_3)$ is a multiple of $(x(P),y(P),z(P))$.
So this can happen at most at one point.
Hence $\dim H^0(Y,\shE \otimes \O_Y(-P)) =1$ and
$H^0(Y, \shE \otimes \O_Y(-P-Q))=0$.
Therefore $\shE \otimes \O_Y(-P)$ is normalized and the $e$-invariant
of $\PP(x,y;z)$ is
$$e= -\deg \, (\shE \otimes \O_Y(-P)) = -3+2=-1 \, .$$
\end{example}

\section{Examples over the complex numbers $\CC$}
\label{complex}

The ruled surfaces together with their forcing sections arising
from tight closure problems yield analytically interesting
examples over the field of complex numbers $K=\CC$.

\begin{corollary}
Let $F$ be a homogeneous polynomial of degree three
such that $R=\CC[x,y,z]/(F)$ defines an
elliptic curve $Y= \Proj\, R$. Let $f,g,h$ be homogeneous
such that $f$ and $g$ are parameters,
$\deg \, h = \deg\, f + \deg\, g$ and $h \not\in (f,g)$.
Then the complement of the forcing section $\PP(f,g)$ in
the ruled surface $\PP(f,g;h)$ over
$Y$ is not affine, but it is a complex Stein space.
The same is true for the open subset
$D(R_+) \subset \Spec\, R[T_1,T_2]/(fT_1+gT_2+h)$.
\end{corollary}

\proof
The condition on the degrees shows that
$h \in (f,g)^\pasoclo$ and that the open complement is not affine.
Since $h \not\in (f,g)$, the forcing sequence, which is an
extension of $\O_Y$ by $\O_Y$, does not split.
Hence due to \cite{umemura} the complement is Stein.
The corresponding statement for the subset in the cone follows.
\qed

\begin{example}
\label{steintwo}
Let $R=\CC[x,y,z]/(x^3+y^3+z^3)$ and $f=x,g=y,h=z^2$.
Then $z^2 \not\in (x,y)$, but $z^2 \in (x,y)^\pasoclo$.
The open subset
$$D(x,y) \subseteq \Spec\, \CC[x,y,z][T_1,T_2]/(x^3+y^3+z^3,\, xT_1+yT_2+z^2)$$
is not affine, but it is a complex Stein space.
\end{example}

\begin{remark}
The first construction of a non-affine but Stein variety
was given by Serre using non-split extensions
of $\O_Y$ by $\O_Y$ on an elliptic curve, see \cite{umemura}
or \cite{bingener}.
Thus we may consider this classical construction
as a construction using forcing algebras.

On the other hand we have to remark that tight closure
takes into account the subtile difference
between affine and Stein.
This shows that tight closure is a conception
of algebraic geometry, not of complex analysis.
\end{remark}

It is also possible to construct new counterexamples
to the hypersection problem. The first counterexample
was given by Coltoiu and Diederich in \cite{coltoiudiederich}.
For this and related problems in complex analysis
see \cite{diederich}. 

\begin{proposition}
Let $R$ be a standard graded normal two-dimensional
$\CC$-algebra, let $f,g$ and $h$ be homogeneous elements in $R$
such that $V(f,g)=V(R_+)$, $h \not\in (f,g)$ and
$\deg f + \deg g - \deg h < 0$.
Then
$$W=D(R_+) \subset \Spec \, R[T_1,T_2]/(fT_1+gT_2+h)$$
is not Stein {\rm(}considered as a complex space{\rm )},
but it fulfills the assumption in the hypersection problem,
i.e. for every analytic surface $S \subset \Spec\, R[T_1,T_2]/(fT_1+gT_2+h)$
the intersection $W \cap S$ is Stein.
\end{proposition}
\proof
Due to \ref{superhoehe} the superheight of $W$ is one,
and \cite[Theorem 5.1]{brennersuperheight} gives that the assumption of the
hypersection problem holds.
The self intersection number of the forcing section on the
corresponding ruled surface is negative,
thus due to \cite{grauertmod} this section is contractible as a complex space
and therefore its complement is not Stein.
Hence the subset $W$ is also not Stein,
because it is a $\CC^\times$-bundle over this complement.
\qed

\begin{example}
Consider
$R=\CC[x,y,z]/(x^4+y^4+z^4)$ and the forcing algebra for the
elements
$x,y;z^3$, hence
$A=\CC[x,y,z,T_1,T_2]/(x^4+y^4+z^4,xT_1+yT_2+z^3)$.
Then $z^3 \not\in (x,y)$ in $R$, but $z^3 \in (x,y)^\pasoclo$,
for the degree of the self intersection is $-4$.
Therefore
$W=D(R_+) \subseteq \Spec \, A$ is not Stein,
but for every analytic surface $S \subset \Spec A$
the intersection $W \cap S$ is Stein.
\end{example}

\section{Plus closure in Positive characteristic}
\label{positive}
The theorem of Smith \cite{smithparameter}, \cite[Theorem 7.1]{hunekeparameter}
says that the tight closure and the plus closure of a
parameter ideal are the same.
In this section we will give a proof for this
in the two-dimensional graded case within our setting.
Let $K$ be an algebraically closed field of positive characteristic $p$
and let $f$ and $g$ be homogeneous parameters in
a two-dimensional standard graded normal $K$-algebra $R$.
Let $h$ be another homogeneous element.
If the complement of the forcing divisor $Z \subset \PP(f,g;h)$
is not affine, then we must find a curve on the ruled surface
disjoined to $Z$.

Suppose first that $ \leadno = \deg\, f + \deg\, g - \deg\, h < 0$.
Then the cohomology class
$\frac{h}{fg}$ has degree $- \leadno $ and becomes
after a Frobenius morphism
$(\frac{h}{fg})^q \in H^1(Y,\O_Y(-q \leadno ))$, but
$H^1(Y,\O_Y(- q \leadno))=0$ for $q=p^{e}$ sufficently large.
Therefore the forcing sequence splits after a Frobenius morphism
and $h$ belongs to the Frobenius closure of $(f,g)$.
Thus we have to consider the case $\leadno = 0$.

\begin{proposition}
\label{frobeniusartinschreier}
Let $K$ be an algebraically closed field of characteristic $p >0$.
Let $R$ be a standard graded normal two-dimensional $K$-algebra
and let $f,g,h$ be homogeneous elements such
that $f$ and $g$ are parameters and such that
$\deg h =\deg f + \deg g$.
Then there exists a composition of a Frobenius morphism and an
Artin-Schreier extension of $Y= \Proj \,R$ such that the
image of the cohomology class $h/fg \in H^1(Y,\O_Y)$
vanishes.
\end{proposition}

\proof
The Frobenius morphism $\Phi$ acts on
$H^1(Y,\O_Y)$ $p$-linear yielding the so called
Fitting decomposition
$H^1(Y,\O_Y) = V_s \oplus V_n$ such that
$\Phi|V_s$ is bijective and $\Phi|V_n$ is nilpotent
\cite[III \S 3]{haramp}.
Thus we may write $c=c_1+c_2$,
where $c_2$ becomes zero after applying a certain power of the Frobenius.
Thus we may assume that $c=c_1 \in V_s$.
Consider the Artin-Schreier sequence
$$0 \lra \ZZ/ (p) \lra \O_Y \stackrel{\Phi - id}{\lra}  \O_Y \lra 0 \, ,$$
which is exact in the $\acute{\rm e}$tale topology. It yields the
exact sequence
$$ 0 \lra H^1_{et}(Y, \ZZ/(p)) \lra H^1(Y,\O_Y)
\stackrel{\Phi - id}{\lra} H^1(Y,\O_Y) \lra \ldots \, \, .$$
There exists a basis $c_j$ of $V_s$ such that $\Phi(c_j)=c_j$,
see \cite[\S 14]{mumfordabelian}.
Thus we may assume that $\Phi(c)-c=0$ and we consider
$c \in H^1_{et}(Y,\ZZ/(p))$.
Hence $c$ represents an Artin-Schreier extension $Y'$ of $Y$ and
the cohomology class $c$ vanishes on $Y'$, see \cite[III. \S 4]{milne}.
\qed

\begin{remark}
\label{artinschreier}
We describe the Artin-Schreier extension appearing in the
last proof explicitly.
Let $c=h/f_1f_2$ and suppose that $c^p-c=0$ in $H^1(Y,\O_Y)$.
This means that
$c^p-c= a_2-a_1$ where
$a_i \in \Gamma(U_i, \O_Y)$, $U_i = D_+(f_i)=\Spec\, R_i$.
Let $U_i'= \Spec\, R_i[T_i]/(T_i^p-T_i+a_i)$, $i=1,2$.
The transition function $T_1 \mapsto T_2+c$ is due to
$$(T_2+c)^p -(T_2+c) +a_1 = T_2^p- T_2 +c^p-c +a_1= T_2^p- T_2+a_2$$
well defined and $U_1'$ and $U_2'$ glue together to a scheme $Y' \ra Y$.
The cohomology class in $Y'$ is $c=h/f_1f_2=T_1-T_2$,
$T_i \in \Gamma(U_i', \O_{Y'})$, so that $c=0$ in $H^1(Y', \O_{Y'})$.
\end{remark}

\begin{remark}
Let $Y = \Proj \, R$ as in {\rm \ref{frobeniusartinschreier}}.
The $p$-rank of $Y$ is
$ \dim V_s \leq g(Y)$.
This is the same as the $p$-rank of the jacobian of $Y$,
see \cite[\S 15]{mumfordabelian}.
The $p$-rank is 0 if and only if the plus closure (=tight closure)
of any homogeneous parameter ideal is
the same as its Frobenius closure.

Suppose now that $Y= \Proj \, K[x,y,z]/(F)$ is an elliptic curve,
thus $H^1(Y, \O_Y) \cong K$.
An elliptic curve with $p$-rank $0$ is
called supersingular (or is said to have Hasse invariant $0$).
The criterion \cite[Proposition IV.4.21]{haralg} says that
$Y$ is supersingular if and only if the
coefficient of $(xyz)^{p-1}$ in $F^{p-1}$ is 0,
or equivalently $F^{p-1} \in (x^p,y^p,z^p)$.

On the other hand,
if the Hasse invariant is 1, then $F^{p-1} \not\in (x^p,y^p,z^p)$,
and the criterion of Fedder \cite[Theorem 3.7]{hunekeparameter}
tells us that
$K[x,y,z]/(F)$ is Frobenius pure.
\end{remark}

\section{Primary Relations and $e$-invariant}
\label{primary}

We give some further estimates of the $e$-invariant of a ruled surface
arising from forcing parameter data which depend upon
the existence of homogeneous primary relations
of some total degree. We will apply this to tight closure problems
of higher rank in the next section.

\begin{lemma}
\label{filtrieren}
Let $R$ be a normal standard graded $K$-algebra.
Let $f_1, \ldots ,f_n \in R$ be homogeneous of degrees $d_i$ and such
that the
$f_i, i \neq j$ are primary for every $j$.
Let $g_i$ be a homogeneous primary relation for $f_i$
of total degree $k$.
Then there exists a sequence
$0 \ra \O_Y \ra \shR(k) \ra \shL \ra 0$ such that
$\shL$ is locally free.
\end{lemma}
\proof
The relation
$(g_1, \ldots ,g_n)$ is a global element $\neq 0$
in
$\Gamma(Y,\shR (k))$ and yields a subsheaf
$\O_Y \subseteq \shR (k)$.
We show that the quotient is locally free, and consider $R_g,\, g=g_1$
(the $D_+(g_i)$ cover $Y$, since they are primary).
$\Gamma(D_+(g), \shR (k))$ is the kernel of the mapping
$$ (R_g)_{e_1} \oplus  \ldots \oplus (R_g)_{e_n}
\stackrel{\sum f_i}{ \lra} (R_g)_k \, $$
($d_i +e_i=k$).
If
$(h_1,\ldots,h_n)$ is an element of the kernel, then
$$(h_1,\ldots,h_n) =
\frac{h_1}{g_1} (g_1, \ldots ,g_n)
+ (0,h_2-\frac{h_1}{g_1} g_2, \ldots,  h_n- \frac{h_1}{g_1} g_n) \, .$$
The second summand is a relation for $(f_2, \ldots, f_n)$.
Because these are also primary, $\shL$ is locally free on $Y$.
\qed

\begin{proposition}
\label{erelationestimate}
Let $R$ be a normal two-dimensional standard graded $K$-al\-gebra.
Let $f_1,f_2,f_3 \in R$ be homogeneous elements
which are pairwise primary of degrees $d_1,d_2,d_3$.
Suppose that $g_1,g_2,g_3 \in R$
is a primary homogeneous relation of total degree $k$.
Set $a= \max(k-d_3,d_1+d_2-k)$.
Then for the normalizing number $\nonu$
and the $e$-invariant of $\PP(f_1,f_2,f_3)$ we have
$$ \nonu \leq a  \deg\, H \, \, 
\mbox{ and }\, \, e \leq |2k-d_1-d_2-d_3| \deg \, H.$$
If $2k=d_1+d_2+d_3$, then
the sheaf of linear forms $\shF(-k)$ is normalized and
the $e$-invariant of $\PP(f_1,f_2,f_3)$ is $0$.
\end{proposition}

\proof
From \ref{filtrieren} we have the sequence
$0 \ra \O_Y \ra \shR (k) \ra \shL \ra 0$,
where $\shL$ is an invertible sheaf.
We know that
$\det \shR(k) \cong \O_Y(-d_1-d_2-d_3+2k)$ from \ref{forcingsheaf3}
and therefore
$\shL =\O_Y(-d_1-d_2-d_3+2k)$.
The dual sequence for the linear forms is then
$0 \ra \O_Y(d_1+d_2+d_3-2k) \ra \shF (-k) \ra \O_Y \ra 0$.
If $d_1+d_2-k \leq k-d_3$, then the sheaf on the left has degree $\leq 0$,
hence $\shF (-k) \otimes \shM$
has no global sections $\neq 0$ whenever
$\shM$ has negative degree.
Hence $\nonu \leq (k-d_3) \deg \, H$.

If $d_1+d_2-k > k-d_3$, we tensorize with $\O_Y(-d_1-d_2-d_3+2k)$
and now $\shF (k- d_1-d_2-d_3)$ is in between two invertible sheaves
of degree $\leq 0$.
Hence $\nonu \leq (d_1+d_2-k) \deg \, H$.
Thus by lemma \ref{lunde} the $e$-invariant is
$e = 2 \nonu - \leadno \deg \, H \leq ( 2a - \leadno) \deg\, H$,
which gives the result.

If $d_1+d_2+d_3=2k$,
then we have the forcing sequence
$0 \ra \O_Y \ra \shF(-k) \ra \O_Y \ra 0$,
thus the sheaf $\shF (-k)$ has sections $\neq 0$
and $\shF(-k)$ is normalized with $\deg \, \shF(-k)=0$.
\qed

\begin{corollary}
\label{bastel}
Let $f_1,f_2,f_3 \in K[x,y,z]$ be homogeneous of degrees $d_1,d_2,d_3$
such that $d_1+d_2+d_3$ is even. Let $g_1,g_2,g_3 \in K[x,y,z]$
be homogeneous of degrees $e_1,e_2,e_3$ and such that
$k=e_i+d_i=(d_1+d_2+d_3)/2$.
Let $F=f_1g_1+f_2g_2+f_3g_3$ and set $R=K[x,y,z]/(F)$.
Suppose that the $f_i$ are pairwise parameters for $R$, that $R$ is normal and
that $V(g_1,g_2,g_3)=V(x,y,z)$.
Then the $e$-invariant of the ruled surface
$\PP(f_1,f_2,f_3)$ over $\Proj\, R$ is $e=0$.
\end{corollary}

\proof
All the conditions in proposition \ref{erelationestimate} are fulfilled.
\qed

\begin{example}
\label{fermat}
Consider a Fermat polynomial $x^m+y^m+z^m \in K[x,y,z]$
and let $R=K[x,y,z]/(x^m+y^m+z^m)$.
Let $f=x^{d_1}, \, g=y^{d_2}, h=z^{d_3}$ such that $d_i < m$.
If $d_1+d_2+d_3 =2m$, then the $e$-invariant
of $\PP(x^{d_1},y^{d_2},z^{d_3})$ is $0$.
Just take $(x^{m-d_1},y^{m-d_2},z^{m-d_3})$ as a primary relation.
\end{example}

\section{Projective bundles of higher rank over a curve}
\label{higherrank}

In this last section we consider again
a two-dimensional standard graded normal $K$-algebra $R$
over an algebraically closed field $K$,
but now we look at the tight closure of three homogeneous primary
elements $(f_1,f_1,f_3)$.
A forth element $f_4$ gives the projective bundle $\PP(f_1,f_2,f_3;f_4)$
of rank two over the smooth base curve $Y= \Proj \, R$ together with
the forcing subbundle $Z=\PP(f_1,f_2,f_3)$, which is itself a ruled surface
over $Y$.
We will need properties of these ruled surfaces
to obtain results on $\PP(f_1,f_2,f_3;f_4)$.
The third self intersection number of the forcing
divisor is $Z^3=(d_1+d_2+d_3-2d_4) \deg\, H$, where $H$
is the hyperplane section on $Y$.

\begin{lemma}
\label{numerischeffektiv}
Let $R$ be a normal standard graded two-dimensional $K$-algebra.
Let $f_1, \ldots ,f_4$ be homogeneous elements
of degrees $d_i$
such that $f_1$ and $f_2$ are parameters.
Set $\leadno =d_1+d_2+d_3-2d_4$.
Suppose that
the $e$-invariant of the forcing subbundle
$Z=\PP(f_1,f_2,f_3)$ is $e \geq 0$.
If $\leadno \deg\, H \geq e $ holds, then $Z$ is a numerically effective
divisor.
If $\leadno \deg\, H > e $, then the pull back $Z|_Z$ is ample.
\end{lemma}
\proof
The intersection of $Z$ with a curve $C \not\subseteq Z$ is $\geq 0$.
The intersection of a curve $ \subset Z$ may be
computed with the pull back of $Z$ on $Z$.
We have $(Z|_Z)^2=Z^3= \leadno \deg \, H$.
Then the result follows from \ref{amplerulecrit}
(The ample criterion of \ref{amplerulecrit} is true even
if $D=C_0 +bF$ is a priori not effective and a similar argument
shows that $D^2 \geq 0$ is equivalent with $D$ numerically effective).
\qed

\begin{corollary}
\label{numerischampel}
Let $K$ be an algebraically closed field and let
$R$ be a normal standard graded two-dimensional $K$-algebra.
Let $f_1, \ldots ,f_4$ be homogeneous elements
such that $f_1$ and $f_2$ are parameters and
such that $\leadno \deg\, H > e \geq 0$,
where $e$ denotes the $e$-invariant of the forcing subbundle
$Z=\PP(f_1,f_2,f_3)$.
Then $f_4 \in (f_1, f_2,f_3)^\pasoclo$ if and only if
$f_4 \in (f_1,f_2,f_3)^{+{\rm gr}}$.
\end{corollary}
\proof
This follows from \ref{pullbackample} and \ref{numerischeffektiv}.
\qed

\begin{theorem}
\label{anwendung}
Let $R$ be a normal standard graded two-dimensional $K$-algebra.
Let $f_1,f_2,f_3$ be homogeneous elements such that
$f_1,f_2$ are parameters and such that the
$e$-invariant of $\PP(f_1,f_2,f_3)$ is $e = 0$.
Let $d_1 +d_2+d_3$ be even and let $m =(d_1 + d_2 +d_3)/2 $.
Then $R_{\geq m} \subseteq (f_1,f_2,f_3)^\pasoclo$.
\end{theorem}
\proof
Let $f_4 $ be homogeneous of degree $m$.
Due to \ref{numerischeffektiv} we know that the forcing divisor
$Z \subset \PP(f_1,f_2,f_3;f_4)$ is numerically effective.
On the other hand the third self intersection number
of $Z$ is zero.
Due to \ref{bignef} the forcing divisor is not big and the
complement is not affine.
\qed

\begin{corollary}
\label{bastel2}
Let $R$ be a normal two-dimensional standard graded $K$-algebra.
Let $f_1,f_2,f_3 \in R$ be homogeneous elements, which are
pairwise primary of degrees $d_1,d_2,d_3$.
Suppose that $d_1+d_2+d_3$ is even and set $m=(d_1+d_2+d_3)/2$.
Suppose that there exists a primary homogeneous relation
of total degree $m$ for $f_1,f_2,f_3$.
Then $R_{\geq m} \subseteq (f_1,f_2,f_3)^\pasoclo$.
\end{corollary}
\proof
This follows from \ref{erelationestimate} and \ref{anwendung}.
\qed

\begin{example}
Let $F=x^m+y^m+z^m$ and $R=K[x,y,z]/(F)$.
Then $R_{\geq m} \subseteq (x^{d_1}, y^{d_2}, z^{d_3})^\pasoclo$,
where $d_1+d_2+d_3=2m$ and $d_i < m$.

For instance we get
$xyz \in (x^2,y^2,z^2)^\pasoclo$ modulo $x^3+y^3+z^3=0$.
This was stated in \cite{hunekeparameter} as an elementary
example of what is not known
in tight closure theory. The first proof was given in \cite{singh}.
\end{example}

In the present case of a projective bundle of rank two
we cannot characterize $f_4 \not\in (f_1,f_2,f_3)^\pasoclo$
by the ampleness of the forcing divisor, as the following
example shows.
The first example of an affine open subset in a three-dimensional
smooth projective variety with no ample divisor
on the complement was given by Zariski and described in \cite{goodman}.

\begin{example}
Let $R=K[x,y]$ and consider on $\PP^1=\Proj\, R$
the projective bundle of rank two
defined by the forcing data
$x^4,y^4,x^4;x^3y^3$.
The third self intersection number of the forcing subbundle $Z$
is zero, hence $Z$ is not ample.
But the complement of $Z$ is affine,
since $x^3y^3 \not\in (x^4,y^4) =(x^4,y^4)^\pasoclo 
= (x^4,y^4,x^4)^\pasoclo $ in the regular ring $K[x,y]$.
Therefore $Z$ is also big.

$Z$ is not numerically effective:
set $e_4=0$ (and $e_1=e_2=e_3=2$), then $Z$ is a hyperplane section.
The pull back of $Z$ on $Z$ yields the hyperplane section
on the ruled surface
$Z= \PP(x^4,y^4,x^4)$ for this grading.
Let $E=Z|_Z$ denote this hyperplane section, let $C =\PP(x^4,y^4) \subset Z$
be the forcing section
and let $L$ be a disjoined section corresponding to $x^4 \in (x^4,y^4)$.
Then we know that $C \sim E + 2\pi^*H$ (\ref{forcingsequence2}(iii))
and therefore
$E.L = C.L - 2\pi^* H.L = - 2\pi^* H.L < 0$.

The divisor $Z$ is also not semiample:
We know that there exists a curve  $L$ such that
$Z.L <0$. Let $P \in L$ and suppose
there exists an effective divisor
$D \sim aZ$ such that $P \not\in D$. Then
$L \not\subset D$ yields a contradiction.
\end{example}

The next results deal with the plus closure in positive
characteristic.

\begin{theorem}
\label{primarynumcrit}
Let $K$ denote an algebraically closed field
of characteristic $p > 0$ and let
$R$ be a normal two-dimensional standard graded $K$-algebra.
Let $f_1,f_2,f_3 \in R$ be homogeneous elements which are
pairwise primary.
Let $g_1,g_2,g_3 \in R$ be a primary homogeneous relation of
total degree $k$. Then for $m= \max\, (k, d_1+d_2+d_3-k)$ we have
$R_{\geq m} \subseteq (f_1,f_2,f_3)^{+{\rm gr}}$.
\end{theorem}
\proof
Let $f_4 \in R_m$, and set $e_4=0$.
The forcing sequence for the relations,
$0 \ra \shR (m) \ra \shR' (m) \ra \O_Y \ra 0$, corresponds to an element
$c \in H^1(Y, \shR (m))$.
We have to show that there exists a smooth projective
curve $\tilde{Y}$ and a finite morphism
$\psi: \tilde{Y} \ra Y$ such that
$\psi^*(c) \in H^1(\tilde{Y}, \O_{\tilde{Y}} )$ is $0$.
The primary relation yields an exact sequence
$0 \ra \O_Y \ra \shR (k) \ra \shL \ra 0$ due to \ref{filtrieren}.
Hence
$$\det \shR (k) = \O_Y(-d_1-d_2-d_3+2k)
\cong \shL \, .$$
Now $ -d_1-d_2-d_3+2k \geq 0$ if and only if
$k \geq (d_1+d_2+d_3-k)$.
If $k< (d_1+d_2+d_3-k)$, we tensorize this sequence with
$\O_Y(d_1+d_2+d_3-2k)$ and get
$0 \ra \O_Y(d_1+d_2+d_3-2k) \ra \shR (d_1+d_2+d_3-k) \ra \O_Y \ra 0 $.

In both cases we have an exact sequence
$0 \ra \O_Y(a) \ra \shR (m) \ra \O_Y(b) \ra 0$
such that $a,b \geq 0$.
Due to \ref{frobeniusartinschreier} and the preceeding remarks there
we know that for a cohomology class $c=\frac{h}{gf} \in H^1(Y, \O_Y(b))$
with $b \geq 0$ there exists a finite mapping
$\varphi : Y' \ra Y$ such that $\varphi^*(c)=0$ in $H^1(Y', \O_{Y'}(b))$.
We apply this first to the image of $c \in H^1(Y, \shR (m))$
in $H^1(Y, \O_Y(b))$ and we may therefore assume
that the image of the cohomology class
$c$ in $H^1(Y', \O_{Y'}(b) )$ vanishes.
Hence we may consider $c \in H^1(Y', \O_{Y'}(a))$ and again
this vanishes after a finite mapping.
\qed

\begin{corollary}
\label{primarynumcor}
Let $K$ be an algebraically closed field of positive characteristic,
let $f_1,f_2,f_3 \in K[x,y,z]$ be homogeneous elements
of degrees $d_1,d_2,d_3$ and let $g_1,g_2,g_3 \in K[x,y,z]$
be primary homogeneous elements of degrees $e_1,e_2,e_3$ such that
$m=e_i+d_i$ and $2m \geq d_1+d_2+d_3$.
Set $F=f_1g_1+f_2g_2+f_3g_3$ and $R=K[x,y,z]/(F)$ and suppose that $R$
is normal and that the $f_i$ are pairwise primary in $R$.
Then $R_{\geq m} \subseteq (f_1,f_2,f_3)^{+{\rm gr}}$.
\end{corollary}
\proof
This follows from \ref{primarynumcrit}, since $m \geq d_1+d_2+d_3-m$.
\qed

\begin{example}
\label{fermatdegree}
Let $F=x^m+y^m+z^m$.
Then $R_{\geq m} \subseteq (x^{d_1}, y^{d_2}, z^{d_3})^{+{\rm gr}}$
for $d_1+d_2+d_3 \leq 2m$ and $d_i < m$.
\end{example}

\begin{example}
Let $F \in K[x,y,z]$ be a homogeneous equation for an elliptic curve.
When is $xyz  \in (x^2,y^2,z^2)^+$ in $R=K[x,y,z]/(F)$?
If the coefficient of $F$ in $xyz$ is not zero, then of course
$xyz \in (x^2,y^2,z^2)$.
Thus we may write $F=Qx^2+Py^2+Sz^2$ so that
$(Q,P,S)$ is a homogeneous relation for $(x^2,y^2,z^2)$.
If $V(Q,P,S)=V(R_+)$, then $xyz \in (x^2,y^2,z^2)^+$.
But what about $F=x^3+y^3+(x+y)z^2$?
\end{example}

\begin{example}
Let $R=K[x,y,z]/(F)$, where $F=x^4+y^4+z^4$.
Then $x^2y^2 \in (x^3,y^3,z^2)^+$.
\end{example}

Suppose now that the normalizing number $\nonu$ of the forcing divisor
$\PP(V) \subset \PP(V')$ is $ >0$.
Then the forcing divisor is big and
there exists a linearly equivalent effective divisor
$D \sim  \PP(V)$ such that its complement is affine, see
\ref{bigaffine}.
The intersection of $\PP(V)$ and $D$ contains a lot of subtile information
for the tight closure problem.

\begin{proposition}
\label{nupositiv}
Let $R$ be a normal two-dimensional standard graded $K$-alge\-bra.
Let $f_1, f_2,f_3$ be homogeneous primary elements
and let $f_0$ be another homogeneous element.
Suppose that $\nonu > 0$.
Then there exists an effective divisor $D$,
$Z \sim D=H+F$, where $H$ is the horizontal component
and $F$ the fiber components.
Moreover, the following hold.

\renewcommand{\labelenumi}{(\roman{enumi})}
\begin{enumerate}

\item 
If $H-Z \cap H$ is not affine, then $f_0 \in (f_1,f_2,f_3)^\pasoclo$.

\item
If $H-Z \cap H$ is affine
{\rm(}this is fulfilled when the pull back $Z|_H$ is ample or
when $H \cap Z$ contains components which lie in a fiber{\rm)},
then there does not exist a finite graded solution
for the tight closure problem, i.e. $f_0 \not\in (f_1,f_2,f_3)^{+{\rm gr}}$.
\end{enumerate}
\end{proposition}

\proof
The condition $\nonu >0$ means that there exists
a positive divisor $L \subset Y$ such that
there exists an effective divisor $D' \sim Z -\pi^*L$.
Then $Z \sim D=D' + \pi^*L$ may be written
as $D= H+F$, where $H$ is a projective subbundle
and $F$ consists of fiber components.

We look at the intersection $Z \cap H$.
(i). $H-Z\cap H = H \cap (\PP(V') - Z) \subseteq \PP(V')-Z$
is a closed subscheme. Hence, if $H- Z \cap H$ is not affine,
then also $\PP(V')-Z$ is not affine and $f_0 \in (f_1,f_2,f_3)^\pasoclo$.

(ii).
We have to show that the forcing divisor
$Z$ intersects every curve $C \not\subseteq Z$ positively.
Due to \ref{effectivevertreter} we only have to consider curves
on $H$. Then the assumption gives the result.
\qed

\begin{example}
Let $Y= \PP^1_K= \Proj\, K[x,y]$ and consider the projective bundles
corresponding to the forcing data $x,y,1;1$.
Then $\nonu >0$ and $Z$ is big, $Z^3=2 >0$ and $Z$ is numerically
effective, but there exists a disjoined curve to $Z$
(the solution section corresponding to $1 \in (x,y,1)$)
and its complement is not affine.

We set $e_4=0$, thus $e_1=e_2=-1$. Eliminating $T_3$ in the forcing
equations yields the splitting forcing sequence
$$0 \lra \AA_Y(1) \times \AA_Y(1) \lra \AA_Y(1) \times \AA_Y(1) \times \AA_Y
\lra \AA_Y \lra 0 \, .$$
The forcing subbundle is $Z \cong \PP^1 \times \PP^1$ given by
the equation $T_4=0$.
We have $\shF'(0) = \O_Y(1) \oplus \O_Y(1) \oplus \O_Y$
and also $\shF'(0) \otimes \O_Y(-1)$ has sections $\neq 0$.
Therefore $\nonu >0$ and $Z$ is big.
A section is for example $xT_1$, thus a divisor $D$
linearly equivalent to $Z$ is given by $D=H+F$, where
$H=\{ T_1 =0 \}$ and $F= \{ x=0 \}$.
$H$ is a Hirzebruch surface $\PP( \AA_Y \times \AA_Y(-1))$
(the blowing up of a projective plane).
$H \cap Z$ is a (horizontal) fiber on $Z \cong \PP^1 \times \PP^1$
and a line on $H$ not meeting the exceptional divisor
(which is also the solution section).
The self intersection number of $Z \cap H$ on $H$ is
$$ (Z|_H)^2=  Z^2.H =Z^2(Z-F)=Z^3-Z^2.F= 2-1=1 \, ,$$
hence $Z$ is numerically effective.
(The self intersection of $Z \cap H$ on $Z$ is 0.)
\end{example}

\begin{example}
Let $K$ be an algebraically closed field of positive characteristic
$p \geq 3$ and consider
$$R=K[x,y,z]/(x^4+ay^4+bz^4+cxz^3+dyz^3) \,  $$
where $a,b,c,d \neq 0 \, $ are chosen such that $Y= \Proj \, R$ is smooth.
We want to show that both cases described
in \ref{nupositiv} do actually occur depending on the coefficients.
$y$ and $z$ are para\-meters and
we consider the elements
$x^4,\, y^4,\, z^4$ and $ xy^2z^3$.
First, let the homogeneous forcing algebra
$$A = R[T_1,T_2,T_3,T_4]/(x^4T_1+ y^4T_2 +z^4T_3 + xy^2z^3T_4)$$
be graded by $e_4=1$ ($e_1=e_2=e_3=3$).
From the curve equation and the homogeneous forcing equation we get
$$ z^3(- (bz+cx+dy)T_1 +zT_3 +xy^2T_4) = y^4(aT_1 -T_2)\, .$$
This gives us the global linear form (of total degree 7) given by
$G$
$$-\frac{bz+cx+dy}{y^4}T_1+\frac{z}{y^4}T_3 + \frac{x}{y^2}T_4 \mbox{ on }
D_+(y)\, \mbox{ and }\,
+\frac{a}{z^3}T_1 - \frac{1}{z^3}T_2 \mbox{ on } D_+(z) \, ,$$
showing that $\nonu > 0$.
We change the grading and set $e_4=0$ and we consider the
linear form $zG$ of total degree 6.
We have then $\PP(V) =V_+(T_4) \sim V_+(G) + V_+(z)=H+F$ on $\PP(V')$
and we are in the situation of \ref{nupositiv}.

When does the intersection $ \PP(V) \cap V_+(z)$ have fiber components?
If $z \neq 0$, then the equation for $G$ on $D_+(y)$ does not vanish,
thus there cannot be fiber components.
So look at $z=0$. Then $x^4+ay^4 =0$ and
the equation for $G$ becomes just $\frac{cx+dy}{y^4}T_1=0$.
Thus there exists a fiber component if and only if
$cx+dy=0=x^4+ay^4$ has a solution, and this means $(d/c)^4=-a$.

\smallskip
Consider the equation $x^4-y^4+z^4+xz^3+yz^3=0$.
This yields a smooth curve and the
intersection has fiber components, hence every curve
intersects the forcing divisor and therefore
$ xy^2z^3 \not\in (x^4,y^4,z^4)^{+{\rm gr}}$.
Does it belong to the tight closure?

The equation $x^4+y^4+z^4+xz^3+yz^3=0$ yields also a smooth curve,
and here the intersection does not have a fiber component.
Hence the intersection is a section and its self intersection number
on $H=V_+(G)$ is negative. Therefore the complement of
it cannot be affine, hence the complement of the forcing
divisor $\PP(V)$ is not affine, thus
$xy^2z^3 \in (x^4,y^4,z^4)^\pasoclo$.
Does it also belong to the (graded) plus closure?
(We have the primary relation $(z,az,bz+cx+cy)$ of total degree $5$,
hence we only know that $R_7 \subset (x^4,y^4,z^4)^{+{\rm gr}}$
due to \ref{primarynumcrit})

\end{example}

%===========================================================

\end{document}